\documentclass[9pt,journal,twocolumn]{IEEEtran}

\usepackage[english]{babel}
\usepackage[latin1]{inputenc}
\usepackage[T1]{fontenc}
\usepackage{url}
\usepackage{graphicx}
\usepackage{subfigure}
\usepackage{color}
\usepackage{amsmath,amsthm,amsfonts,amssymb}
\usepackage{colortbl}
\usepackage{url}
\usepackage{booktabs}
\usepackage{epsfig}
\usepackage{pstricks}
\usepackage{pst-node}
\usepackage{pst-grad}
\usepackage{ifthen}
\usepackage{algorithm}
\usepackage{algpseudocode}
\usepackage{float}
\newfloat{algorithm}{t}{lop}

\newboolean{draft}
\newcommand{\isdraft}[2]{\ifthenelse{\boolean{draft}}{#1}{#2}}

\isdraft{\usepackage{setspace}}{}                              % DRAFT
\isdraft{\usepackage[footnotesize]{caption}}{}                 % DRAFT
\isdraft{\usepackage{paralist}}{}                              % DRAFT

\theoremstyle{definition}

\theoremstyle{plain}
\newtheorem{Theorem}{Theorem}
\newtheorem{Corollary}{Corollary}

\newtheorem*{Assumptions}{Assumptions}

\definecolor{red}{rgb}{0.69921875,0.28125,0.28125}
\definecolor{green}{rgb}{0.265625,0.609375,0.265625}
\definecolor{blue}{rgb}{0.3828125,0.3828125,0.69140625}
\definecolor{lightred}{rgb}{0.9140625,0.8046875,0.8046875}
\definecolor{lightblue}{rgb}{0.87109375,0.87109375,0.9296875}
\definecolor{lightgreen}{rgb}{0.796875,0.91796875,0.796875}
\definecolor{lightyellow}{RGB}{232,238,184}
\definecolor{white}{RGB}{253,253,253}

\newcommand{\mypar}[1]{{\bf #1.}}

\isdraft{\renewcommand{\baselinestretch}{1.38}}{}

\title{D-ADMM: A Communication-Efficient Distributed Algorithm For Separable Optimization}

\author{Jo\~ao F.~C.~Mota, Jo\~ao M.~F.~Xavier, Pedro M.~Q.~Aguiar, and~Markus P\"uschel
\IEEEcompsocitemizethanks{\IEEEcompsocthanksitem Jo\~ao F.~C.~Mota, Jo\~ao, M.~F.~Xavier, and Pedro M.~Q.~Aguiar  are with Instituto de Sistemas e Rob\'otica (ISR), Instituto Superior T\'ecnico (IST), Technical University of Lisbon, Portugal.
\IEEEcompsocthanksitem Markus P{\"u}schel is with the Department of
Computer Science at ETH Zurich, Switzerland.%
\IEEEcompsocthanksitem Jo\~ao F.~C.~Mota is also with the Department of Electrical and Computer Engineering at Carnegie Mellon University, USA.}
\thanks{This work was supported by the following grants from Funda\c{c}\~{a}o para a Ci\^{e}ncia e Tecnologia (FCT): CMU-PT/SIA/0026/2009, PEst-OE/EEI/LA0009/2011, and SFRH/BD/33520/2008 (through the Carnegie Mellon/Portugal Program managed by ICTI).}}

\begin{document}

\maketitle

\begin{abstract}
  We propose a distributed algorithm, named Distributed Alternating Direction Method of Multipliers (D-ADMM), for solving separable optimization problems in networks of interconnected nodes or agents. In a separable optimization problem there is a private cost function and a private constraint set at each node. The goal is to minimize the sum of all the cost functions, constraining the solution to be in the intersection of all the constraint sets. D-ADMM is proven to converge when the network is bipartite or when all the functions are strongly convex, although in practice, convergence is observed even when these conditions are not met. We use D-ADMM to solve the following problems from signal processing and control: average consensus, compressed sensing, and support vector machines. Our simulations show that D-ADMM requires less communications than state-of-the-art algorithms to achieve a given accuracy level. Algorithms with low communication requirements are important, for example, in sensor networks, where sensors are typically battery-operated and communicating is the most energy consuming operation.
\end{abstract}


\isdraft{\pagebreak}{}

\section{Introduction}
\label{Sec:intro}

	In this paper, we propose a distributed algorithm for solving \textit{separable optimization problems}:
	\begin{equation}\label{Eq:IntroSeparableProb}
    \begin{array}{ll}
      \underset{x}{\text{minimize}} & f_1(x) + f_2(x) + \cdots + f_P(x) \\
      \text{subject to} & x \in X_1 \cap X_2 \cap \cdots \cap X_P\,,
    \end{array}
  \end{equation}
	where~$x \in \mathbb{R}^n$ is the global optimization variable, and~$x^\star$ will denote any solution of~\eqref{Eq:IntroSeparableProb}. As illustrated in Fig.~\ref{Fig:IntroFig7Nodes}, we associate a network of~$P$ nodes with problem~\eqref{Eq:IntroSeparableProb}, where only node~$p$ has access to its private cost function~$f_p$ and to its private set~$X_p$. Each node can only communicate with its neighbors, but all of them have to solve~\eqref{Eq:IntroSeparableProb} in a cooperative way. We call any method that solves~\eqref{Eq:IntroSeparableProb} without using a central node and without aggregating data at any specific location a \textit{distributed algorithm}.

  \mypar{Contributions} The goal of this paper is twofold: to show that the recent distributed algorithm proposed in~\cite{Mota12-DistributedBP} for a specific problem called Basis Pursuit can be generalized to solve the class~\eqref{Eq:IntroSeparableProb}; and to show that, for many problems of interest, the resulting algorithm requires usually significant less communications than prior distributed algorithms to achieve a given solution accuracy. This also includes algorithms that were specifically designed for a particular problem and are not applicable to the entire problem class~\eqref{Eq:IntroSeparableProb}. Algorithms with low communication cost are relevant, for example, in sensor networks where communication is often the most energy-consuming task and the nodes rely on batteries~\cite{Akyildiz02-WirelessNetworksASurvey,Fischione11-DesignPrinciplesWSN}.

	\isdraft{                                                                  % DRAFT
    \begin{figure}                                                           % DRAFT
      \vspace{-0.8cm}                                                        % DRAFT
      \centering                                                             % DRAFT
      \includegraphics[scale=0.9]{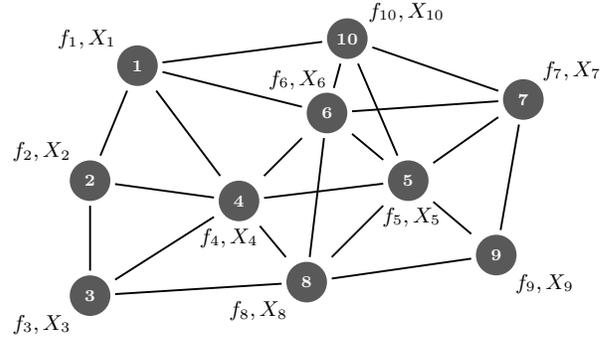}            % DRAFT
      \caption{Network with~$P = 10$ nodes. Node~$p$ only knows~$f_p$ and~$X_p$, but cooperates with its neighbors in order to solve~\eqref{Eq:IntroSeparableProb}.}            % DRAFT
      \label{Fig:IntroFig7Nodes}                                             % DRAFT
      \vspace{-0.6cm}                                                        % DRAFT
    \end{figure}                                                             % DRAFT
  }{
	\begin{figure}
  \centering
  \psscalebox{0.9}{
        \begin{pspicture}(7.0,4.5)
            \def\nodesimp{
							\pscircle*[linecolor=black!65!white](0,0){0.3}
						}

            \rput(1.0,3.7){\rnode{N1}{\nodesimp}}
						\rput(0.3,2.0){\rnode{N2}{\nodesimp}}
						\rput(0.3,0.3){\rnode{N3}{\nodesimp}}
						\rput(2.5,1.7){\rnode{N4}{\nodesimp}}
						\rput(5.0,2.0){\rnode{N5}{\nodesimp}}
						\rput(3.8,3.0){\rnode{N6}{\nodesimp}}
						\rput(6.7,3.2){\rnode{N7}{\nodesimp}}
						\rput(3.5,0.5){\rnode{N8}{\nodesimp}}
						\rput(6.3,0.9){\rnode{N9}{\nodesimp}}
						\rput(4.1,4.1){\rnode{N10}{\nodesimp}}

            \rput(1,3.7){\scriptsize \textcolor{black!3!white}{$\mathbf{1}$}}
						\rput[rb](0.68765248,3.94987802){$f_1,X_1$}
						\rput(0.3,2){\scriptsize \textcolor{black!3!white}{$\mathbf{2}$}}
						\rput[rb](0.017157,2.282843){$f_2,X_2$}
						\rput(0.3,0.3){\scriptsize \textcolor{black!3!white}{$\mathbf{3}$}}
						\rput[rt](0.017157,0.017157){$f_3,X_3$}
						\rput(2.5,1.7){\scriptsize \textcolor{black!3!white}{$\mathbf{4}$}}
						\rput[t](2.35,1.3){$f_4,X_4$}
						\rput(5,2){\scriptsize \textcolor{black!3!white}{$\mathbf{5}$}}
						\rput[t](5.05,1.6){$f_5,X_5$}
						\rput(3.8,3){\scriptsize \textcolor{black!3!white}{$\mathbf{6}$}}
						\rput[rb](3.8,3.35){$f_6,X_6$}
						\rput(6.7,3.2){\scriptsize \textcolor{black!3!white}{$\mathbf{7}$}}
						\rput[lb](6.982843,3.482843){$f_7,X_7$}
						\rput(3.5,0.5){\scriptsize \textcolor{black!3!white}{$\mathbf{8}$}}
						\rput[tr](3.21715729,0.21715729){$f_8,X_8$}
						\rput(6.3,0.9){\scriptsize \textcolor{black!3!white}{$\mathbf{9}$}}
						\rput[lt](6.58284271,0.61715729){$f_9,X_9$}
						\rput(4.1,4.1){\scriptsize \textcolor{black!3!white}{$\mathbf{10}$}}
						\rput[lb](4.41234752,4.34987802){$f_{10},X_{10}$}

						\ncline[nodesep=0.37cm]{-}{N1}{N2}
						\ncline[nodesep=0.37cm]{-}{N1}{N4}
						\ncline[nodesep=0.37cm]{-}{N2}{N3}
						\ncline[nodesep=0.37cm]{-}{N2}{N4}
						\ncline[nodesep=0.37cm]{-}{N4}{N5}
						\ncline[nodesep=0.37cm]{-}{N4}{N6}
						\ncline[nodesep=0.37cm]{-}{N5}{N6}
						\ncline[nodesep=0.37cm]{-}{N5}{N7}
						\ncline[nodesep=0.37cm]{-}{N3}{N4}
						\ncline[nodesep=0.37cm]{-}{N1}{N6}
						\ncline[nodesep=0.37cm]{-}{N6}{N7}
						\ncline[nodesep=0.37cm]{-}{N8}{N4}
						\ncline[nodesep=0.37cm]{-}{N8}{N5}
						\ncline[nodesep=0.37cm]{-}{N8}{N9}
						\ncline[nodesep=0.37cm]{-}{N8}{N3}
						\ncline[nodesep=0.37cm]{-}{N8}{N6}
						\ncline[nodesep=0.37cm]{-}{N9}{N5}
						\ncline[nodesep=0.37cm]{-}{N9}{N7}
						\ncline[nodesep=0.37cm]{-}{N10}{N6}
						\ncline[nodesep=0.37cm]{-}{N10}{N7}
						\ncline[nodesep=0.37cm]{-}{N10}{N1}
						\ncline[nodesep=0.37cm]{-}{N10}{N5}

            %\psgrid
            \end{pspicture}
  }
  \bigskip
  \caption{Network with~$P = 10$ nodes. Node~$p$ only knows~$f_p$ and~$X_p$, but cooperates with its neighbors in order to solve~\eqref{Eq:IntroSeparableProb}.}
  \label{Fig:IntroFig7Nodes}
  \end{figure}
  }

	\mypar{Formal problem statement}
	Given a network with~$P$ nodes, we associate each~$f_p$ and~$X_p$ in~\eqref{Eq:IntroSeparableProb} with the $p$th node of the network. We make the following assumptions:
	\isdraft{                                                                % DRAFT
      \begin{inparaenum}                                                   % DRAFT
        \item Each~$f_p: \mathbb{R}^n \xrightarrow{} \mathbb{R}$ is a convex function over~$\mathbb{R}^n$, and each set~$X_p$ is closed and convex;                               % DRAFT
        \label{Ass:FunctionsAndSets}                                       % DRAFT
        \item Problem~\eqref{Eq:IntroSeparableProb} is solvable;           % DRAFT
        \label{Ass:ProblemSolvable}                                        % DRAFT
        \item The network is connected and it does not vary with time;     % DRAFT
        \label{Ass:Network}                                                % DRAFT
        \item A coloring scheme of the network is available.               % DRAFT
        \label{Ass:Coloring}                                               % DRAFT
       \end{inparaenum}                                                    % DRAFT
  }{
	\begin{Assumptions}
	\hfill

		\begin{enumerate}
			\item Each~$f_p: \mathbb{R}^n \xrightarrow{} \mathbb{R}$ is a convex function over~$\mathbb{R}^n$, and each set~$X_p$ is closed and convex.
			\label{Ass:FunctionsAndSets}

			\item Problem~\eqref{Eq:IntroSeparableProb} is solvable.
			\label{Ass:ProblemSolvable}

			\item The network is connected and it does not vary with time.
			\label{Ass:Network}

			\item A coloring scheme of the network is available.
			\label{Ass:Coloring}
		\end{enumerate}
  \end{Assumptions}
  }

  Assumption~\ref{Ass:ProblemSolvable}) implies that~\eqref{Eq:IntroSeparableProb} has at least one solution~$x^\star$. In Assumption~\ref{Ass:Network}), a network is connected if there is a path between every pair of nodes. Finally, in Assumption~\ref{Ass:Coloring}), a coloring scheme is an assignment of numbers to the nodes of the network such that no adjacent nodes have the same number. These numbers are usually called colors, and they will be used to set up our distributed algorithm. Note that, in wireless scenarios, coloring schemes are often used in \textit{Media Access Control} (MAC) protocols to determine the nodes' order of communication.

  Under the previous assumptions, we solve the following problem: \textit{given a network, design a distributed algorithm that solves~\eqref{Eq:IntroSeparableProb}}. By ``distributed'' we mean there is no notion of a central or special node and each node communicates only with its neighbors; also, only node~$p$ has access to~$f_p$ or~$X_p$ at any time during or before the algorithm.

  Our solution for this problem relies on the \textit{Alternating Direction Method of Multipliers} (ADMM), which has become very popular in recent years; see~\cite{Boyd11-ADMM} for a survey. Specifically, we use an extended version of ADMM, whose proof of convergence was recently established in~\cite{Han12-NoteOnADMM}. This result will also guarantee the convergence of our algorithm for some problems of interest.

	\mypar{Related work}
	Gradient and subgradient methods, including incremental versions, are long known to yield distributed  algorithms (in the sense defined before); see, e.g., \cite{Tsitsiklis86-DistrAsynchronous,Rabbat04-DistributedOptimizationSensorNetworks,CXVOptimInSPandComm}. Advantages of these methods are computational simplicity at each node and theoretical robustness guarantees. However, they generally require too many iterations (and hence communications) to converge.

	Augmented Lagrangian methods have also been used for distributed optimization, e.g., \cite{Bertsekas97-ParallelDistributed,Jakovetic11-CooperativeConvexOptimization,Rabbat05-GeneralizedConsensusAlgs}. They consist of two loops: an outer loop updating the dual variables, and an inner loop updating the primal variables. The most common outer loop algorithm is the gradient method, yielding the method of multipliers. For the inner loop, common choices are Gauss-Seidel and Jacobi methods.

	The Alternating Direction Method of Multipliers (ADMM)~\cite{Boyd11-ADMM} is an augmented Lagrangian-based algorithm that consists of only one loop. ADMM is not directly applicable to~\eqref{Eq:IntroSeparableProb}: one has to reformulate that problem first. Possible reformulations were addressed in~\cite{Schizas08-ConsensusAdHocWSNsPartI} and~\cite{Zhu09-DistributedInNetworkChannelCoding}, yielding algorithms that require two and one communication steps per ADMM iteration, respectively. Other work that explores these algorithms for particular instances of~\eqref{Eq:IntroSeparableProb} include~\cite{Bazerque10-DistributedSpectrumSensing,Forero10-ConsensusBasedDistributedSVMs,Mateos10-DistributedSparseLinearRegression,Bazerque11-GroupLasso,Erseghe11-FastConsensusByADMM}. The algorithm we propose is also based on ADMM (on an extended version), but applied to a different reformulation of~\eqref{Eq:IntroSeparableProb}. Our simulations show that the proposed algorithm requires less communications than any of the previous approaches.

	All the above algorithms solve~\eqref{Eq:IntroSeparableProb} in a distributed way. There are, however, other algorithms that solve~\eqref{Eq:IntroSeparableProb}, but are not distributed in our sense. For example, the algorithm in~\cite[\S7.2]{Boyd11-ADMM} solves~\eqref{Eq:IntroSeparableProb}, but it requires an all-to-all communication in each iteration; this can only be accomplished in networks that are fully connected or that have a central node. In contrast, our algorithm and the ones described above are distributed and can run on any connected network topology.

\section{Algorithm Derivation}
\label{Sec:AlgDerivation}
	To derive the algorithm, we reformulate~\eqref{Eq:IntroSeparableProb} to make ADMM applicable. As mentioned before, several reformulations are possible: ours takes advantage of node coloring. First we introduce some notation.

  \mypar{Network notation} Networks are represented as undirected graphs~$\mathcal{G}=(\mathcal{V},\mathcal{E})$, where~$\mathcal{V} = \{1,2,\ldots,P\}$ is the set of nodes and~$\mathcal{E} \subseteq \mathcal{V} \times \mathcal{V}$ is the set of edges. The cardinality of these sets is represented respectively by~$P$ and~$E$. An edge is represented by~$(i,j)$, with~$i<j$, and $(i,j) \in \mathcal{E}$ means that nodes~$i$ and~$j$ can exchange data with each other. We define the neighborhood~$\mathcal{N}_p$ of a node~$p$ as the set of nodes connected to node~$p$, but excluding it; the cardinality of this set, $D_p := |\mathcal{N}_p|$, is the degree of node~$p$.

  \mypar{Coloring} We assume the network is given together with a coloring scheme of~$C$ colors. The set of nodes that have color~$c$ will be denoted with~$\mathcal{C}_c$, for~$c=1,\ldots,C$, and its cardinality with $C_c =|\mathcal{C}_c|$. Note that~$\{\mathcal{C}_c\}_{c=1}^C$ partitions~$\mathcal{V}$.

  \mypar{Problem manipulations}  Without loss of generality, assume the nodes are ordered such that the first~$C_1$ nodes have color~$1$, the next~$C_2$ nodes have color~$2$, and so on, i.e., $\mathcal{C}_1 = \{1,2,\ldots,C_1\}$, $\mathcal{C}_2 = \{C_1+1,C_1+2,\ldots,C_1+C_2\}$, \ldots. We decouple problem~\eqref{Eq:IntroSeparableProb} by assigning copies of the global variable~$x$ to each node and then constrain all copies to be equal. Let~$x_p \in \mathbb{R}^n$ denote the copy held by node~$p$. As in~\cite{Zhu09-DistributedInNetworkChannelCoding}, we constrain all copies to be equal in an edge-based way, and rewrite~\eqref{Eq:IntroSeparableProb} as
%   \begin{equation}\label{Eq:ClonedProb}
%     \begin{array}{cl}
%       \isdraft{\text{minimize}}{\underset{\bar{x} = (x_1,\ldots,x_P)}{\text{minimize}}} & f_1(x_1) + f_2(x_2) +\cdots + f_P(x_P) \\
%       \isdraft{
%         \text{subject to} & x_p \in X_p\,,\, p = 1,\ldots,P \,,\qquad x_i = x_j \,,\, (i,j) \in \mathcal{E}\,,
%       }{
%          \text{subject to} & x_p \in X_p\,,\quad p = 1,\ldots,P \\
%                            & x_i = x_j \,,\quad (i,j) \in \mathcal{E}\,,
%       }
%     \end{array}
%   \end{equation}
  \begin{equation}\label{Eq:ClonedProb}
    \begin{array}{cl}
      \underset{\bar{x} = (x_1,\ldots,x_P)}{\text{minimize}} & f_1(x_1) + f_2(x_2) +\cdots + f_P(x_P) \\
      \isdraft{
        \text{subject to} & x_p \in X_p\,,\, p = 1,\ldots,P \,,\qquad x_i = x_j \,,\, (i,j) \in \mathcal{E}\,,
      }{
         \text{subject to} & x_p \in X_p\,,\quad p = 1,\ldots,P \\
                           & x_i = x_j \,,\quad (i,j) \in \mathcal{E}\,,
      }
    \end{array}
  \end{equation}
  where~$\bar{x}= (x_1,\ldots,x_P) \in (\mathbb{R}^n)^P$ is the optimization variable. Problem~\eqref{Eq:ClonedProb} is no longer coupled by a global variable, as~\eqref{Eq:IntroSeparableProb}, but instead by the new equations~$x_i=x_j$, for all the pairs~$(i,j) \in \mathcal{E}$. These equations enforce all copies to be equal since the network is connected (cf. Assumption~\ref{Ass:Network})). Note that these constraints can be written more compactly as~$(B^\top \otimes I_n) \bar{x} = 0$, where~$B \in \mathbb{R}^{P\times E}$ is the node arc-incidence matrix of the graph, $I_n$ is the identity matrix in~$\mathbb{R}^n$, and $\otimes$ is the Kronecker product. Each column of~$B$ is associated with an edge~$(i,j) \in \mathcal{E}$ and has $1$ and~$-1$ in the $i$th and $j$th entry, respectively; the remaining entries are zeros. Our numbering assumption induces a natural partition of~$B$ as $\begin{bmatrix}B_1^\top & B_2^\top & \cdots & B_C^\top\end{bmatrix}^\top$\!\!, where the columns of~$B_c^\top$ are associated to the nodes
with color~$c$. We partition~$\bar{x}$ similarly: $\bar{x} = (\bar{x}_1,\ldots,\bar{x}_C)$, where~$\bar{x}_c \in (\mathbb{R}^n)^{C_c}$ collects the copies of all nodes with color~$c$. This enables rewriting~\eqref{Eq:ClonedProb} as
% 	\begin{equation}\label{Eq:ClonedProb2}
% 			\begin{array}{ll}
% 				\isdraft{\text{minimize}}{\underset{\bar{x}_1,\ldots,\bar{x}_C}{\text{minimize}}} & \sum_{c=1}^C \sum_{p \in \mathcal{C}_c} f_p(x_p) \\
% 				\isdraft{
%           \text{subject to} & \bar{x}_c \in \bar{X}_c\,,\, c = 1,\ldots,C\,, \qquad \sum_{c=1}^C (B_c^\top \otimes I_n)\bar{x}_c= 0\,,
%         }{
% 				\text{subject to} & \bar{x}_c \in \bar{X}_c\,,\quad c = 1,\ldots,C \\
% 				                  & \sum_{c=1}^C (B_c^\top \otimes I_n)\bar{x}_c= 0\,,
%         }
% 			\end{array}
% 	\end{equation}
	\begin{equation}\label{Eq:ClonedProb2}
      \begin{array}{ll}
        \underset{\bar{x}_1,\ldots,\bar{x}_C}{\text{minimize}} & \sum_{c=1}^C \sum_{p \in \mathcal{C}_c} f_p(x_p) \\
        \isdraft{
          \text{subject to} & \bar{x}_c \in \bar{X}_c\,,\, c = 1,\ldots,C\,, \qquad \sum_{c=1}^C (B_c^\top \otimes I_n)\bar{x}_c= 0\,,
        }{
        \text{subject to} & \bar{x}_c \in \bar{X}_c\,,\quad c = 1,\ldots,C \\
                          & \sum_{c=1}^C (B_c^\top \otimes I_n)\bar{x}_c= 0\,,
        }
      \end{array}
  \end{equation}
  where~$\bar{X}_c := \prod_{p \in \mathcal{C}_c} X_p$. Problem~\eqref{Eq:ClonedProb2} can be solved with the Extended ADMM, explained next.

  \mypar{Extended ADMM}
  The Extended ADMM is a natural generalization of the \textit{Alternating Direction Method of Multipliers} (ADMM) \cite{Han12-NoteOnADMM}. Given~$C$ functions~$g_c$, $C$ sets~$X_c$, and~$C$ matrices~$A_c$, all with the same number of rows, the extended ADMM solves
%   \begin{equation}\label{Eq:ProblemSolvedByExtADMM}
%     \begin{array}{cl}
%       \isdraft{\text{minimize}}{\underset{x_1,\ldots,x_C}{\text{minimize}}} & \sum_{c=1}^C g_c(x_c) \\
%       \isdraft{
%         \text{subject to} & x_c \in X_c\,,\, c = 1, \ldots,C \,, \qquad \sum_{c=1}^C A_c x_c = 0\,,
%       }{
%       \text{subject to} & x_c \in X_c\,,\quad  c = 1, \ldots,C \\
%                         & \sum_{c=1}^C A_c x_c = 0\,,
%       }
%     \end{array}
%   \end{equation}
  \begin{equation}\label{Eq:ProblemSolvedByExtADMM}
    \begin{array}{cl}
      \underset{x_1,\ldots,x_C}{\text{minimize}} & \sum_{c=1}^C g_c(x_c) \\
      \isdraft{
        \text{subject to} & x_c \in X_c\,,\, c = 1, \ldots,C \,, \qquad \sum_{c=1}^C A_c x_c = 0\,,
      }{
      \text{subject to} & x_c \in X_c\,,\quad  c = 1, \ldots,C \\
                        & \sum_{c=1}^C A_c x_c = 0\,,
      }
    \end{array}
  \end{equation}
  where~$x:=(x_1,\ldots,x_C)$ is the optimization variable. The extended ADMM consists of iterating on~$k$:
  \begin{align}
        x_1^{k+1} &= \underset{x_1 \in X_1}{\arg\min} \,\,\, L_{\rho}(x_1, x_2^{k},\ldots,x_P^k;\lambda^k)
				\label{Eq:ExtADMMAlg1}
				\\
        x_2^{k+1} &= \underset{x_2 \in X_2}{\arg\min} \,\,\, L_{\rho}(x_1^{k+1}, x_2, x_3^k, \ldots, x_C^k; \lambda^k)
				\label{Eq:ExtADMMAlg2}
				\\
				&\,\,\,\vdots
				\notag
				\\
				x_C^{k+1} &= \underset{x_C \in X_C}{\arg\min} \,\,\, L_{\rho}(x_1^{k+1}, x_2^{k+1}, \ldots, x_{C-1}^{k+1}, x_C; \lambda^k)
				\label{Eq:ExtADMMAlg3}
				\isdraft{\displaybreak[3]}{}
				\\
        \lambda^{k+1} &= \lambda^k + \rho \sum_{c = 1}^C A_c x_c^{k+1}\,,
        \label{Eq:extADMMAlg4}
  \end{align}
  where
	$%\begin{equation}\label{Eq:AugmLagrangian}
			L_\rho(x;\lambda) = \sum_{c=1}^C \bigl(g_c(x_c) + \lambda^\top A_c x_c\bigr) + \frac{\rho}{2}\bigl\|\sum_{c=1}^C A_cx_c\bigr\|^2\,\,
	$%\end{equation}
  is the augmented Lagrangian of~\eqref{Eq:ProblemSolvedByExtADMM}, $\lambda$ is the dual variable, and~$\rho$ is a positive parameter. When~$C = 2$, \eqref{Eq:ExtADMMAlg1}-\eqref{Eq:extADMMAlg4} becomes the ordinary ADMM and it converges under very mild assumptions. When~$C > 2$, there is only a known proof of convergence when all the functions~$g_c$ are strongly convex~\cite{Han12-NoteOnADMM}. In particular, the following theorem holds.

  \begin{Theorem}[\cite{Mota11-ADMMProof,Han12-NoteOnADMM}]\label{Teo:ConvergenceADMM}
		Let~$g_c:\mathbb{R}^{n_c}\xrightarrow{} \mathbb{R}$ be a convex function over~$\mathbb{R}^{n_c}$, $X_c \subseteq \mathbb{R}^{n_c}$ a closed convex set, and~$A_c$ an $m\times n_c$ matrix, for~$c=1,\ldots,C$. Assume~\eqref{Eq:ProblemSolvedByExtADMM} is solvable and that either
		\isdraft{
       1) $C=2$ and each~$A_c$ has full column-rank, or 2) $C\geq2$ and each~$g_c$ is strongly convex.
    }{
		\begin{enumerate}
      \item $C=2$ and each~$A_c$ has full column-rank,
      \item or $C\geq2$ and each~$f_c$ is strongly convex.
      %, and $\rho < \underset{c=1,\ldots,C}{\min} \{\frac{2\mu_c}{3(C-1)\|A_c\|^2}\}$.
    \end{enumerate}
    }
    Then, the sequence~$\{(x_1^k,\ldots,x_C^k,\lambda^k)\}$ generated by~\eqref{Eq:ExtADMMAlg1}-\eqref{Eq:extADMMAlg4} converges to $(x_1^\star, \ldots, x_C^\star, \lambda^\star)$, where $(x_1^\star, \ldots, x_C^\star)$ solves~\eqref{Eq:ProblemSolvedByExtADMM} and~$\lambda^\star$ solves the dual problem of~\eqref{Eq:ProblemSolvedByExtADMM}:
		$%\begin{equation}\label{Eq:ADMMDualProb}
      \max_\lambda G_1(\lambda) + \cdots + G_C(\lambda)
    $,
    where~$G_c(\lambda) = \inf_{x_c \in X_c} (g_c(x_c) + \lambda^\top A_cx_c)$, for $c=1,\ldots,C$.
  \end{Theorem}
  A proof for case~$1)$ can be found in~\cite{Mota11-ADMMProof}, which generalizes the proofs of~\cite{Bertsekas97-ParallelDistributed,Boyd11-ADMM}. A proof for case~$2)$ can be found in~\cite{Han12-NoteOnADMM}. It is believed that~\eqref{Eq:ExtADMMAlg1}-\eqref{Eq:extADMMAlg4} still converges when~$C>2$ and  each~$A_c$ has full column-rank, i.e., that the generalization of Theorem~\ref{Teo:ConvergenceADMM} under case~$1)$ still holds~\cite{Mota12-DistributedBP,Han12-NoteOnADMM,He12-ADMGaussianBack}. Recently, \cite{Luo12-LinearConvergenceADMM} proved that if we replace~$\rho$ in~\eqref{Eq:extADMMAlg4} by a small constant, the resulting algorithm converges linearly.
  %\begin{Conjecture}[\cite{Han12-NoteOnADMM,He12-ADMGaussianBack,Mota12-DistributedBP}]\label{Conj:ExtendedADMM}
	%	The conclusions of Theorem~\ref{Teo:ConvergenceADMM} still hold for any~$C\geq 2$	whenever problem~\eqref{Eq:ProblemSolvedByExtADMM} is solvable, and for each~$c=1,\ldots,C$, $g_c$ is convex over~$\mathbb{R}^{n_c}$, $X_c$ is closed and convex, and~$A_c$ is a full column-rank matrix.
  %\end{Conjecture}

	\mypar{Applying the extended ADMM}
	We now apply the extended ADMM to problem~\eqref{Eq:ClonedProb2}, which has the format of~\eqref{Eq:ProblemSolvedByExtADMM}. We start by showing that the $c$th optimization problem in~\eqref{Eq:ExtADMMAlg1}-\eqref{Eq:ExtADMMAlg3} yields~$C_c$ optimization problems that can be solved in parallel. For example, $\bar{x}_1 = (x_1,\ldots,x_{C_1})$ is updated as
	\begin{equation}\label{Eq:MinStep1}
		\bar{x}_1^{k+1}
		=
		\underset{\bar{x}_1 \in \bar{X}_1}{\arg\min}\, \sum_{p \in \mathcal{C}_1} f_p(x_p)
		+
		{\lambda^k}^\top \!\!A_1 \bar{x}_1 + \frac{\rho}{2}\Bigl\|A_1 \bar{x}_1 + \sum_{c=2}^C A_c \bar{x}_c^k\Bigr\|^2,
	\end{equation}
	where~$A_1 = B_1^\top \otimes I_n$. The last term in~\eqref{Eq:MinStep1} can be written as
	\begin{equation}\label{Eq:MinStep2}
		\frac{\rho}{2}\bar{x}_1^\top A_1^\top A_1 \bar{x}_1 + \rho\,\bar{x}_1^\top \sum_{c=2}^C A_1^\top A_c \bar{x}_c^k + \frac{\rho}{2}\Bigl\|\sum_{c=2}^C A_c \bar{x}_c^k\Bigr\|^2\,.
	\end{equation}
	In the first term, $A_1^\top A_1 = B_1 B_1^\top \otimes I_n$, where $B_1 B_1^\top$ is a diagonal block of the graph Laplacian. Since the nodes with color~$1$ are not neighbors between themselves, $B_1 B_1^\top$ will be a diagonal matrix, with the degrees of the respective nodes in the diagonal. This means $\bar{x}_1^\top A_1^\top A_1 \bar{x}_1 = \sum_{p \in \mathcal{C}_1} D_p \|x_p\|^2$. Similarly, in the second term, $A_1^\top A_c = B_1 B_c^\top \otimes I_n$, where $B_1 B_c^\top$ corresponds to an off-diagonal block of the Laplacian matrix. For $i\neq j$, the $ij$th entry of the Laplacian matrix contains~$-1$ if nodes~$i$ and~$j$ are neighbors, and~$0$ otherwise. This implies $\bar{x}_1^\top \sum_{c=2}^C A_1^\top A_c \bar{x}_c^k = -\sum_{p \in \mathcal{C}_1} \sum_{j \in \mathcal{N}_p} x_p^\top x_j^k$. Finally, the last term of~\eqref{Eq:MinStep2} does not depend on~$\bar{x}_1$ and can be ignored from the optimization problem. Thus, \eqref{Eq:MinStep1} simplifies to
% 	\isdraft{
%     \begin{equation}\label{Eq:MinStep3}
%     \bar{x}_1^{k+1} = \underset{\bar{x}_1 = (x_1,\ldots,x_{C_1})}{\arg\min}\,\, \sum_{p \in \mathcal{C}_1} \biggl[ f_p(x_p) + \Bigl(\gamma_p^k - \rho \!\!\sum_{j \in \mathcal{N}_p} x_j^k \Bigr)^\top x_p + \frac{\rho D_p}{2} \|x_p\|^2\biggr],
%   \end{equation}
%   }{
% 	\begin{multline}\label{Eq:MinStep3}
%     \bar{x}_1^{k+1} = \underset{\bar{x}_1 = (x_1,\ldots,x_{C_1})}{\arg\min}\,\, \sum_{p \in \mathcal{C}_1} \biggl[ f_p(x_p) + \Bigl(\gamma_p^k - \rho \!\!\sum_{j \in \mathcal{N}_p} x_j^k \Bigr)^\top x_p \\+ \frac{\rho D_p}{2} \|x_p\|^2\biggr],
% 	\end{multline}
%   }
  \begin{equation}\label{Eq:MinStep3}
    \bar{x}_1^{k+1} = \underset{\bar{x}_1 \in \bar{X}_1}{\arg\min}\, \sum_{p \in \mathcal{C}_1} f_p(x_p) + \Bigl(\gamma_p^k - \rho \!\!\sum_{j \in \mathcal{N}_p} x_j^k \Bigr)^\top x_p \\+ \frac{\rho D_p}{2} \|x_p\|^2,
	\end{equation}
	where~$\gamma_p^k := \sum_{j \in \mathcal{N}_p} \lambda_{pj}^k$ comes from the second term in~\eqref{Eq:MinStep1}: $((B_1 \otimes I_n) \lambda^k)^\top \bar{x}_1 = \sum_{p \in \mathcal{C}_1} \sum_{j \in \mathcal{N}_p} {\lambda_{pj}^k}^\top x_p$. We decomposed~$\lambda$ edge-wise: $\lambda = (\ldots,\lambda_{ij},\ldots)$, where~$\lambda_{ij}$ is defined for~$i<j$ and associated to the constraint~$x_i = x_j$ in~\eqref{Eq:ClonedProb}. It is now clear that~\eqref{Eq:MinStep3} decomposes into~$C_1$ problems that can be solved in parallel. For the other colors, we can apply a similar reasoning, but we must be careful defining~$\gamma_p^k$, due to the nodes' relative numbering. Its general definition is~$\gamma_p^k := \sum_{j \in \mathcal{N}_p} \text{sign}(j-p) \lambda_{pj}^k$, where~$\text{sign}(a) = 1$, if~$a\geq 0$, and $\text{sign}(a) = -1$, otherwise. Note that we extended the definition of~$\lambda_{ij}$ for $i>j$ such that~$\lambda_{ij} := \lambda_{ji}$. Algorithm~\ref{Alg:DADMMGeneral} shows the
resulting algorithm, named \textit{Distributed-ADMM}, or D-ADMM.

	\isdraft{                                                                      % DRAFT
	  \vspace{-0.3cm}                                                              % DRAFT
    \begin{minipage}{\textwidth}                                                 % DRAFT
    \begin{minipage}[c]{0.5\linewidth}                                           % DRAFT
    \begin{algorithm}[H]                                                         % DRAFT
    \caption{\small D-ADMM}                                                      % DRAFT
    \algrenewcommand\algorithmicrequire{\textbf{Initialization:}}                % DRAFT
    \label{Alg:DADMMGeneral}                                                     % DRAFT
    \begin{algorithmic}[1]                                                       % DRAFT
    \scriptsize                                                                % DRAFT
    \Require for all~$p \in \mathcal{V}$, set $\gamma_{p}^{1} = x_p^{1} = 0$ and $k=1$
    \Repeat                                                                      % DRAFT
    \For{$c =1,\ldots,C$}                                                        % DRAFT
        \ForAll{$p \in \mathcal{C}_c$ [in parallel]}                             % DRAFT
            $$                                                                   % DRAFT
                v_p^{k} = \gamma_p^{k}-                                          % DRAFT
                \rho \sum_{\begin{subarray}{c}                                   % DRAFT
                             j \in \mathcal{N}_p \\                              % DRAFT
                             j < p                                               % DRAFT
                           \end{subarray}                                        % DRAFT
                }x_j^{k+1} - \rho \sum_{\begin{subarray}{c}                      % DRAFT
                             j \in \mathcal{N}_p \\                              % DRAFT
                             j > p                                               % DRAFT
                           \end{subarray}                                        % DRAFT
                }x_j^{k}                                                         % DRAFT
            $$                                                                   % DRAFT
            \vspace{-0.3cm}                                                      % DRAFT
            \label{SubAlg:GeneralDADMMv}                                         % DRAFT
        \State and find                                                          % DRAFT
            $$                                                                   % DRAFT
            x_p^{k+1} = \begin{array}[t]{cl}                                     % DRAFT
                            \textrm{argmin} & f_p(x_p) + {v_p^{k}}^\top x_p + \frac{D_p \rho}{2}\|x_p\|^2\\
                            \textrm{s.t.} & x_p \in X_p
                          \end{array}
           $$                                                                    % DRAFT
           \label{SubAlg:ADMMGenProb}                                            % DRAFT
        \State Send~$x_p^{k+1}$ to $\mathcal{N}_p$                               % DRAFT
        \label{SubAlg:ADMMGenComm}                                               % DRAFT
    \EndFor                                                                      % DRAFT
    \EndFor                                                                      % DRAFT
                                                                                 % DRAFT
    \ForAll{$p \in \mathcal{V}$ [in parallel]} \vspace{0.15cm}                   % DRAFT
    \hfill                                                                       % DRAFT

        $                                                                        % DRAFT
            \gamma_p^{k+1} = \gamma_p^{k} + \rho \sum_{j \in \mathcal{N}_p} (x_p^{k+1} -  x_j^{k+1})
        $\label{SubAlg:ADMMGenDualVarUpdt}  \vspace{0.15cm}                      % DRAFT
                                                                                 % DRAFT
    \EndFor                                                                      % DRAFT
    \State $k \gets k+1$                                                         % DRAFT
    \Until{some stopping criterion is met}                                       % DRAFT
    \end{algorithmic}                                                            % DRAFT
    \end{algorithm}                                                              % DRAFT
    \vspace{-0.3cm}                                                              % DRAFT
    \end{minipage}                                                               % DRAFT
    \hfill                                                                       % DRAFT
    \begin{minipage}[c]{0.4\linewidth}                                           % DRAFT
       %\begin{figure}                                                           % DRAFT
        \centering                                                               % DRAFT
        \includegraphics[scale=0.9]{figures/SingleCol/MatrixPartition.eps}       % DRAFT
        \vspace{0.3cm}                                                           % DRAFT
        \captionof{figure}{\footnotesize Row partition and column partition of~$A$ into~$P$ blocks. A block in the row (resp. column) partition is a set of rows (resp. columns).}                     % DRAFT
        \label{Fig:PartitionOfA}                                                % DRAFT
       %\end{figure}                                                             % DRAFT
    \end{minipage}                                                               % DRAFT
    \end{minipage}                                                               % DRAFT
  }{
	\begin{algorithm}[H]
    \caption{D-ADMM}
    \algrenewcommand\algorithmicrequire{\textbf{Initialization:}}
    \label{Alg:DADMMGeneral}
    \begin{algorithmic}[1]
    \small
    \Require for all~$p \in \mathcal{V}$, set $\gamma_{p}^{1} = x_p^{1} = 0$ and $k=1$
    \Repeat
    \For{$c =1,\ldots,C$}
        \ForAll{$p \in \mathcal{C}_c$ [in parallel]}
            $$
                v_p^{k} = \gamma_p^{k}-
                \rho \sum_{\begin{subarray}{c}
                             j \in \mathcal{N}_p \\
                             j < p
                           \end{subarray}
                }x_j^{k+1} - \rho \sum_{\begin{subarray}{c}
                             j \in \mathcal{N}_p \\
                             j > p
                           \end{subarray}
                }x_j^{k}
            $$
            \vspace{-0.3cm}
            \label{SubAlg:GeneralDADMMv}
        \State and find
            $$
            x_p^{k+1} = \begin{array}[t]{cl}
                            \textrm{argmin} & f_p(x_p) + {v_p^{k}}^\top x_p + \frac{D_p \rho}{2}\|x_p\|^2\\
                            \textrm{s.t.} & x_p \in X_p
                          \end{array}
           $$
           \label{SubAlg:ADMMGenProb}
        \State Send~$x_p^{k+1}$ to $\mathcal{N}_p$
        \label{SubAlg:ADMMGenComm}
    \EndFor
    \EndFor

    \ForAll{$p \in \mathcal{V}$ [in parallel]} \vspace{0.15cm}
    \hfill

        $
            \gamma_p^{k+1} = \gamma_p^{k} + \rho \sum_{j \in \mathcal{N}_p} (x_p^{k+1} -  x_j^{k+1})
        $\label{SubAlg:ADMMGenDualVarUpdt}  \vspace{0.15cm}

    \EndFor
    \State $k \gets k+1$
    \Until{some stopping criterion is met}
    \end{algorithmic}
  \end{algorithm}
  }
	In Algorithm~\ref{Alg:DADMMGeneral}, the edge-wise dual variables~$\lambda_{ij}$ were totally replaced by the node-wise dual variables~$\gamma_p$. This is because the problem in step~\ref{SubAlg:ADMMGenProb} depends only on~$\gamma_p^k$ and not on the individual $\lambda_{ij}^k$'s. The update for~$\gamma_p$ in step~\ref{SubAlg:ADMMGenDualVarUpdt} stems from replacing $\lambda_{ij}^{k+1} = \lambda_{ij}^k + \text{sign}(j-i)(x_i^{k+1} - x_j^{k+1})$ in the definition of~$\gamma_p^{k+1}$.

	Algorithm~\ref{Alg:DADMMGeneral} is asynchronous in the sense that nodes operate in a color-based order, with nodes with the same color operating in parallel. Since nodes with the same color are not neighbors, we would apparently need some kind of coordination to execute the algorithm. Actually, such coordination is not needed provided each node knows its own color and the colors of its neighbors. In fact, as soon as node~$p$ has received~$x_j^{k+1}$ from all its neighbors with lower colors, node~$p$ can ``work,'' since step~\ref{SubAlg:ADMMGenProb} (and subsequently step~\ref{SubAlg:ADMMGenComm}) can be performed. In conclusion, knowing its own and its neighbors' colors provides an automatic coordination mechanism.	Regarding the convergence of D-ADMM, we have:
	\begin{Corollary}\label{Cor:DADMMConv}
		Let Assumptions~\ref{Ass:FunctionsAndSets})\! -\! \ref{Ass:Coloring}) hold. Then, Algorithm~\ref{Alg:DADMMGeneral} produces a sequence~$(x_1^k,\ldots,x_P^k)$ convergent to~$(x^\star, \ldots,x^\star)$, where~$x^\star$ solves~\eqref{Eq:IntroSeparableProb}, whenever $1$) the network is bipartite, or $2$) each~$f_p$ is strongly convex.
		%\begin{enumerate}
		%	\item the network is bipartite;
		%	\item each~$f_p$ is strongly convex;
		%	\item Conjecture~\ref{Conj:ExtendedADMM} is true.
		%\end{enumerate}
	\end{Corollary}
	\begin{proof}
		The proof is based on showing that the conditions of Theorem~\ref{Teo:ConvergenceADMM} are satisfied. First, note that Assumptions~\ref{Ass:FunctionsAndSets}) and~\ref{Ass:ProblemSolvable}) and the equivalence between~\eqref{Eq:IntroSeparableProb} and~\eqref{Eq:ClonedProb2} imply that problem~\eqref{Eq:ClonedProb2} is solvable, that each function $\sum_{p \in \mathcal{C}_c} f_p(x_p)$ is convex over~$\mathbb{R}^n$, and that each set~$\bar{X}_c$ is closed and convex. Now, we will see that Assumption~\ref{Ass:Network}) implies that each~$B_c^\top \otimes I_n$ has full column-rank. Note that it is sufficient to prove that~$B_c^\top$ has full column-rank. If, on the other hand, we prove that~$B_cB_c^\top$ has full rank, then the result follows because~$\text{rank}(B_cB_c^\top) = \text{rank}(B_c^\top)$. Note that $B_cB_c^\top$ is a diagonal matrix, where the diagonal contains the degrees of the nodes belonging to the subnetwork composed by the nodes in~$\mathcal{C}_c$. Since no node has degree~$0$ (cf. Assumption~\ref{Ass:Network})), $B_cB_c^\top$ has full rank.

		Finally, note that a bipartite network can be colored with just two colors. In that case, condition~$1$) of Theorem~\ref{Teo:ConvergenceADMM} is satisfied together with the remaining conditions, which ensures the convergence of Algorithm~\ref{Alg:DADMMGeneral}. When the network is non-bipartite and each~$f_p$ is strongly convex, we are in case~$2$) of Theorem~\ref{Teo:ConvergenceADMM}, which again ensures the convergence of Algorithm~\ref{Alg:DADMMGeneral}.
		%If Conjecture~\ref{Conj:ExtendedADMM} is true, we just proved that all of its conditions are satisfied (in particular that each~$B_c^\top \otimes I_n$ has full column-rank).
	\end{proof}

\section{Applications}
\label{Sec:Applications}

	We will now see how some important optimization problems can be recast as~\eqref{Eq:IntroSeparableProb}. These reformulations, except the one for LASSO, are not new: see~\cite{Zhu09-DistributedInNetworkChannelCoding,Bazerque10-DistributedSpectrumSensing,Mateos10-DistributedSparseLinearRegression,Forero10-ConsensusBasedDistributedSVMs}. Therefore, we refer to these references for the details of solving the optimization problem in step~\ref{SubAlg:ADMMGenProb} of Algorithm~\ref{Alg:DADMMGeneral}. We note that in all the problems, except in consensus, none of the functions~$f_p$ is strongly convex. Therefore, D-ADMM is only guaranteed  to converge under condition $1$) of Corollary~\ref{Cor:DADMMConv}. Nevertheless, in section~\ref{Sec:Experiments}, we will see that in practice D-ADMM, not only converges for all these problems, but also outperforms previous work in terms of the number of communications, including the ADMM-based algorithms~\cite{Zhu09-DistributedInNetworkChannelCoding,Schizas08-ConsensusAdHocWSNsPartI,Mateos10-DistributedSparseLinearRegression}.

	\mypar{Consensus} Consensus is a fundamental problem in networks~\cite{Erseghe11-FastConsensusByADMM,Oreshkin10-OptimizationAnalysisDistrAveraging}. Given a network with~$P$ nodes, node~$p$ generates a number, say~$\theta_p$, and the goal is to compute the average~$\theta^\star = (1/P)\sum_{p=1}^P \theta_p$ at every node. Consensus can be cast as~\cite{Rabbat04-DistributedOptimizationSensorNetworks,Erseghe11-FastConsensusByADMM}:
	\isdraft{                                                                     % DRAFT
  $                                                                             % DRAFT
      \text{minimize} \,\,\, \frac{1}{2}\sum_{p=1}^P (x - \theta_p)^2\,,        % DRAFT
  $                                                                             % DRAFT
  }{
	$$
      \underset{x}{\text{minimize}} \,\,\, \frac{1}{2}\sum_{p=1}^P (x - \theta_p)^2\,,
  $$
  }
  which is clearly an unconstrained version of~\eqref{Eq:IntroSeparableProb}, with~$f_p(x) = (1/2)(x-\theta_p)^2$; thus, it can be solved with D-ADMM. In this case, the problem of step~\ref{SubAlg:ADMMGenProb} of Algorithm~\ref{Alg:DADMMGeneral} has a closed-form solution: $x_p^{k+1} = (\theta_p - v_p^k)/(1+D_p\rho)$.

  \mypar{Sparse solutions of linear systems} Finding sparse solutions of linear systems is important in many areas, including statistics, compressed sensing, and cognitive radio~\cite{Tibshirani96-RegressionShrinkageLasso,Bazerque10-DistributedSpectrumSensing}. A common approach to tackle this problem is by solving LASSO~\cite{Tibshirani96-RegressionShrinkageLasso} or BPDN~\cite{AtomicDecompBP}, respectively,
  \isdraft{
		\begin{equation}\label{Eq:LassoBPDN}
			\text{LASSO}\,:\,\,\,\,
			\begin{array}[t]{ll}
			\underset{x}{\text{minimize}} & \|x\|_1 \\
			\text{subject to} & \|Ax - b\| \leq \sigma
			\end{array}
			\qquad
			\text{BPDN}\,:\,\,\,\,
			\underset{x}{\text{minimize}} \,\,\, \|Ax - b\|^2 + \beta\|x\|_1\,,
		\end{equation}
	}{
		\begin{align}
		\text{LASSO:} & \phantom{aaaaa}
				\begin{array}[t]{ll}
						\underset{x}{\text{minimize}} & \|x\|_1 \\
						\text{subject to} & \|Ax - b\| \leq \sigma\,,
				\end{array}
				\label{Eq:Lasso}
	  \\
	  \text{BPDN:} &  \phantom{aaaaa}
			\begin{array}[t]{ll}
				\underset{x}{\text{minimize}} & \|Ax - b\|^2 + \beta\|x\|_1\,,
			\end{array}
    \label{Eq:BPDN}
	\end{align}
	}
	where the matrix~$A \in \mathbb{R}^{m \times n}$, the vector~$b \in \mathbb{R}^m$, and the parameters~$\sigma, \beta>0$ are given. LASSO first appeared in~\cite{Tibshirani96-RegressionShrinkageLasso} to denote a related problem, although the problem in~\isdraft{\eqref{Eq:LassoBPDN}}{\eqref{Eq:Lasso}} is known by the same name. We solve LASSO and BPDN in two different scenarios, visualized in Fig.~\ref{Fig:PartitionOfA}: \textit{row partition} (resp. \textit{column partition}), where each node stores a block of rows (resp. columns) of~$A$. While in the row partition vector~$b$ is partitioned similarly to~$A$, in the column partition we assume all nodes know the full vector~$b$.

	We propose solving LASSO with a column partition and BPDN with a row partition. The reverse cases, i.e., LASSO with a row partition and BPDN with a column partition, cannot be trivially recast as~\eqref{Eq:IntroSeparableProb}. However, in our previous work~\cite{Mota12-DistributedBP}, we solved Basis Pursuit (i.e., LASSO with~$\sigma = 0$) for both the row and the column partition.

	 % Figure: Partition of M
   % ======================
    %\isdraft{                                                                   % DRAFT
    %  \begin{figure}                                                            % DRAFT
    %    \centering                                                              % DRAFT
    %    \includegraphics[scale=0.94]{figures/FigsSingleCol/PartitionOfA.eps}    % DRAFT
    %    %\vspace{-0.3cm}                                                         % DRAFT
    %    \caption{Row partition and column partition of~$A$ into~$P$ blocks. In the row partition a block is a set of rows; in the column partition a block is a set of columns.}      % DRAFT
    %    \label{Fig:PartitionOfA}                                                % DRAFT
    %  \end{figure}                                                              % DRAFT
    %}{
    \isdraft{}{
    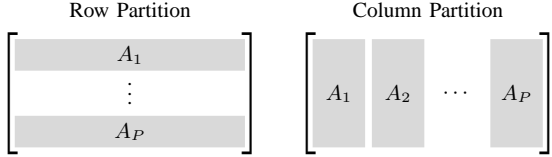
\begin{figure}
    \centering
    \psscalebox{0.9}{
    \begin{pspicture}(8.0,2.5)
				\def\brackets{
					\rput(0,0){\psscalebox{1.5 1.9}{$\Biggl[$}}
				}

				\def\blockmatrix{
					\psframe*[linecolor=black!15!white,fillstyle=solid](0,0)(3.4,0.4667)
				}
				\rput[bl](0.1,0.2){\blockmatrix}
				\rput[bl](0.1,1.334){\blockmatrix}

				\rput(1.8,1.5667){$A_1$}
				\rput(1.8,1.1){$\vdots$}
				\rput(1.8,0.4334){$A_P$}

				\rput(0.01,1){\brackets}
				\rput(3.59,1){\psscalebox{-1 1}{\brackets}}

				\def\blockmatrix{
					\psframe*[linecolor=black!15!white,fillstyle=solid](0,0)(0.775,1.6)
				}

				\rput[bl](4.5,0.2){\blockmatrix}
				\rput[bl](5.3750,0.2){\blockmatrix}
				\rput[bl](7.1250,0.2){\blockmatrix}

				\rput(4.8875,1){$A_1$}
				\rput(5.7625,1){$A_2$}
				\rput(6.6375,1){$\cdots$}
				\rput(7.5125,1){$A_P$}

				\rput(4.41,1){\brackets}
				\rput(7.99,1){\psscalebox{-1 1}{\brackets}}

				\rput(1.8,2.23){Row Partition}
				\rput(6.2,2.23){Column Partition}

				%\psgrid
    \end{pspicture}
		}
		\vspace{-0.1cm}
    \caption{Row partition and column partition of~$A$ into~$P$ blocks. A block in the row (resp. column) partition is a set of rows (resp. columns).}
    \label{Fig:PartitionOfA}
    \end{figure}
    }

    \mypar{LASSO: column partition} Assume~$A$ is partitioned by columns, and the $p$th block is only known at node~$p$. Also, assume vector~$b$, parameter~$\sigma$, and the number of nodes~$P$ are available at all nodes. LASSO in this scenario cannot be directly recast as~\eqref{Eq:IntroSeparableProb}: we will have to do it through duality. However, only solving the ordinary dual of LASSO will not allow us to recover a primal solution afterwards, since its objective is not strictly convex. We thus start by regularizing LASSO, making it strictly convex:
%     \begin{equation}\label{Eq:LassoRegularized}
%       \begin{array}{ll}
%         \isdraft{\text{minimize}}{\underset{x}{\text{minimize}}} & \|x\|_1 + \frac{\delta}{2} \|x\|^2 \\
%         \text{subject to} & \|Ax - b\| \leq \sigma\,,
%       \end{array}
%     \end{equation}
    \begin{equation}\label{Eq:LassoRegularized}
      \begin{array}{ll}
        \underset{x}{\text{minimize}} & \|x\|_1 + \frac{\delta}{2} \|x\|^2 \\
        \text{subject to} & \|Ax - b\| \leq \sigma\,,
      \end{array}
    \end{equation}
    where~$\delta >0$ is small enough. This regularization is inspired by~\cite{Friedlander07ExactRegularizationConvexPrograms}, which establishes exact regularization conditions. By exact, we mean there exists~$\bar{\delta} > 0$ such that the solution of~\eqref{Eq:LassoRegularized} is always a LASSO solution, for~$\delta \leq \bar{\delta}$. One of these conditions is that the objective is linear and the constraint set is the intersection of a linear system with a closed polyhedral cone. Although LASSO can be recast as
	\isdraft{                                                             % DRAFT
    \begin{equation}\label{Eq:LassoManipulated}                         % DRAFT
    \begin{array}{ll}                                                   % DRAFT
      \underset{x,t,u,v}{\text{minimize}} & 1_n^\top t \\               % DRAFT
      \text{subject to} & \|u\| \leq v\,,\qquad  u = Ax - b\,,\qquad v = \sigma\,,\qquad  -t \leq x \leq t\,,
    \end{array}                                                         % DRAFT
    \end{equation}                                                      % DRAFT
  }{
	\begin{equation}\label{Eq:LassoManipulated}
    \begin{array}{ll}
      \underset{x,t,u,v}{\text{minimize}} & 1_n^\top t \\
      \text{subject to} & \|u\| \leq v \\
                        & u = Ax - b\,,\,\,\, v = \sigma \\
                        & x \leq t\,,\,\,\, -x\leq t\,,
    \end{array}
  \end{equation}
  }
  where~$1_n \in \mathbb{R}^n$ is the vector of ones, the closed convex cone $\{(u,v)\,:\, \|u\|\leq v\}$ is not polyhedral; thus, there is not a proof of exact regularization for~\eqref{Eq:LassoManipulated}. However, experimental results in~\cite{Friedlander07ExactRegularizationConvexPrograms} suggest that exact regularization might occur for non-polyhedral cones. In the simulations discussed in the next section, we solved~\eqref{Eq:LassoRegularized} always with $\delta=10^{-2}$ and the corresponding solutions never differed more than $0.5\%$ from the ``true'' solution of LASSO.

  We now introduce a variable~$y \in \mathbb{R}^{m}$ in~\eqref{Eq:LassoRegularized} and rewrite it as:
  \isdraft{
    \begin{equation}\label{Eq:LassoWithY}
    \begin{array}{ll}
      \underset{x,y}{\text{minimize}} & \sum_{p=1}^P (\|x_p\|_1 + \frac{\delta}{2}\|x_p\|^2) \\
      \text{subject to} & \|y\| \leq \sigma \,,\qquad  y = \sum_{p=1}^P A_p x_p - b\,.
    \end{array}
  \end{equation}
  }{
  \begin{equation}\label{Eq:LassoWithY}
    \begin{array}{ll}
      \underset{x,y}{\text{minimize}} & \sum_{p=1}^P (\|x_p\|_1 + \frac{\delta}{2}\|x_p\|^2) \\
      \text{subject to} & \|y\| \leq \sigma \\
                        & y = \sum_{p=1}^P A_p x_p - b\,.
    \end{array}
  \end{equation}
  }
  If we only dualize the last constraint of~\eqref{Eq:LassoWithY}, we get the dual problem of minimizing $\sum_{p=1}^P g_p(\lambda)$, where~$g_p(\lambda) := \frac{1}{P}(b^\top \lambda + \sigma\|\lambda\|) - \inf_{x_p}(\|x_p\|_1 + (A_p^\top \lambda)^\top x_p + \frac{\delta}{2}\|x_p\|^2)$ is the function associated to node~$p$. This problem is clearly an unconstrained version of~\eqref{Eq:IntroSeparableProb}.

  %We point out that problem~\eqref{Eq:LassoRegularized} can be ill-conditioned due to the small value that has to be used for~$\delta$, if we want to get a good approximation of the original problem. In particular, the problem solved at each node can be very time-consuming, since it may require lots of internal iterations. We used~$\delta = 10^{-3}$ in our simulations. While such value for~$\delta$ required a considerable amount of iterations for solving the internal problem of each node (we set~$500$ as the number of maximum iterations), it also allowed to get solutions of the original problem~\eqref{Eq:Lasso} with relative errors of~$0.1\%$, using few communications. BPDN, addressed next, handles the same type of problems as LASSO, but its distributed implementation does not require solving ill-conditioned problems.

  \mypar{BPDN: row partition} In BPDN, $A$ and~$b$ are partitioned by rows, with the blocks~$A_p$ and~$b_p$ stored at node~$p$. In this scenario, BPDN can be readily rewritten as
  \begin{equation}\label{Eq:BPDN2}
		\underset{x}{\text{minimize}} \,\,\, \sum_{p=1}^P \Bigl( \|A_px - b_p\|^2 + \frac{\beta}{P}\|x\|_1\Bigr)\,,
  \end{equation}
  which is an unconstrained version of~\eqref{Eq:IntroSeparableProb}: just make~$f_p(x) := \|A_px - b_p\|^2 + \frac{\beta}{P}\|x\|_1$.
  %The corresponding problem in step~\ref{SubAlg:ADMMGenProb} of Algorithm~\ref{Alg:DADMMGeneral} can be solved efficiently, for example, with the GPSR solver~\cite{GPSR}.

  \mypar{Distributed support vector machines} A \textit{Support Vector Machine} is an optimization problem that arises in machine learning in the context of classification and regression~\cite[Ch.7]{Bishop06-PatternRecognitionMachineLearning}. While there are several possible formulations for an SVM, here we solve~\cite[\S7.1.1]{Bishop06-PatternRecognitionMachineLearning}
  \isdraft{
    \begin{equation}\label{Eq:SVM}
    \begin{array}{ll}
      \underset{s,r,\xi}{\text{minimize}} & \frac{1}{2}\|s\|^2 + \beta \,1_K^\top \xi \\
      \text{subject to} & y_k \, (s^\top x_k - r) \geq 1 - \xi_k\,,\,k = 1,\ldots,K \,,\qquad \xi \geq 0\,,
    \end{array}
  \end{equation}
  }{
  \begin{equation}\label{Eq:SVM}
		\begin{array}{ll}
			\underset{s,r,\xi}{\text{minimize}} & \frac{1}{2}\|s\|^2 + \beta \,1_K^\top \xi \\
			\text{subject to} & y_k \, (s^\top x_k - r) \geq 1 - \xi_k\,,\quad k = 1,\ldots,K \\
												& \xi \geq 0\,,
		\end{array}
  \end{equation}
  }
  where the parameter~$\beta > 0$ and the pairs~$(x_k,y_k)$, $k=1,\ldots,K$, are given. Each point~$x_k$ belongs to one of two classes: $y_k = 1$ or~$y_k = -1$. The goal in solving~\eqref{Eq:SVM} is to find an hyperplane~$\{x \in \mathbb{R}^n\,:\, s^\top x = r\}$ that best separates the two classes. The optimization variables in~\eqref{Eq:SVM} are~$s \in \mathbb{R}^n$, the vector orthogonal to the hyperplane, $r \in \mathbb{R}$, the hyperplane offset, and~$\xi \in \mathbb{R}^K$, the vector of slack variables. We assume that~$K$ is divisible by the number of nodes~$P$, and that each node knows~$m := K/P$ pairs of points~$(x_k,y_k)$. The resulting problem can be formulated as an unconstrained version of~\eqref{Eq:IntroSeparableProb} by setting
  $$
		f_p(s,r) =
		\begin{array}[t]{cl}
			\underset{\bar{\xi}_p}{\inf} & \frac{1}{2P}\|s\|^2 + \beta \,1_{m}^\top \bar{\xi}_p \\
			\text{s.t.} & Y_p \, (X_p s - r 1_m) \geq 1_m - \bar{\xi}_p \\
			& \bar{\xi}_p \geq 0\,,
		\end{array}
	$$
	where~$Y_p$ is a diagonal matrix with the $y_k$'s corresponding to node~$p$ in the diagonal, and $X_p$ is an $m\times n$ matrix with each row containing~$x_k^\top$. Note that the size of the variable to be transmitted among the nodes is~$n+1$, corresponding to the size of the global variable~$(s,r)$: the variables~$\bar{\xi}_p$ are internal to each node.

\section{Simulation Results}
\label{Sec:Experiments}

	This section shows simulation results of D-ADMM and related algorithms solving the problems presented in the previous section. We focus on the ADMM-based algorithms~\cite{Schizas08-ConsensusAdHocWSNsPartI} and~\cite{Zhu09-DistributedInNetworkChannelCoding}, since these are the best among the distributed algorithms for~\eqref{Eq:IntroSeparableProb}. We start by discussing the performance measure.

	\isdraft{                                                                      % DRAFT
    \begin{figure}                                                               % DRAFT
      \vspace{-0.2cm}
      \centering                                                                 % DRAFT
      \includegraphics[scale=0.95]{figures/SingleCol/ExperimentalResults.eps}    % DRAFT
      \subfigure{\label{SubFig:Consensus}}                                       % DRAFT
      \subfigure{\label{SubFig:Lasso}}                                           % DRAFT
      \subfigure{\label{SubFig:BPDN}}                                            % DRAFT
      \subfigure{\label{SubFig:SVM}}                                             % DRAFT
      \vspace{-0.7cm}                                                            % DRAFT
      \caption{Results of the simulations for $\text{(a)}$ consensus, $\text{(b)}$ LASSO, $\text{(c)}$ BPDN, and $\text{(d)}$ SVM.}                                                         % DRAFT
     \label{Fig:Experiments}                                                     % DRAFT
     %\vspace{-0.9cm}                                                             % DRAFT
    \end{figure}                                                                 % DRAFT
  }{
	\begin{figure*}
     \centering
     \subfigure[Consensus]{\label{SubFig:Consensus}
     \begin{pspicture}(0.23\linewidth,5.1cm)
       \rput(0.117\linewidth,2.6){\includegraphics[scale=0.24]{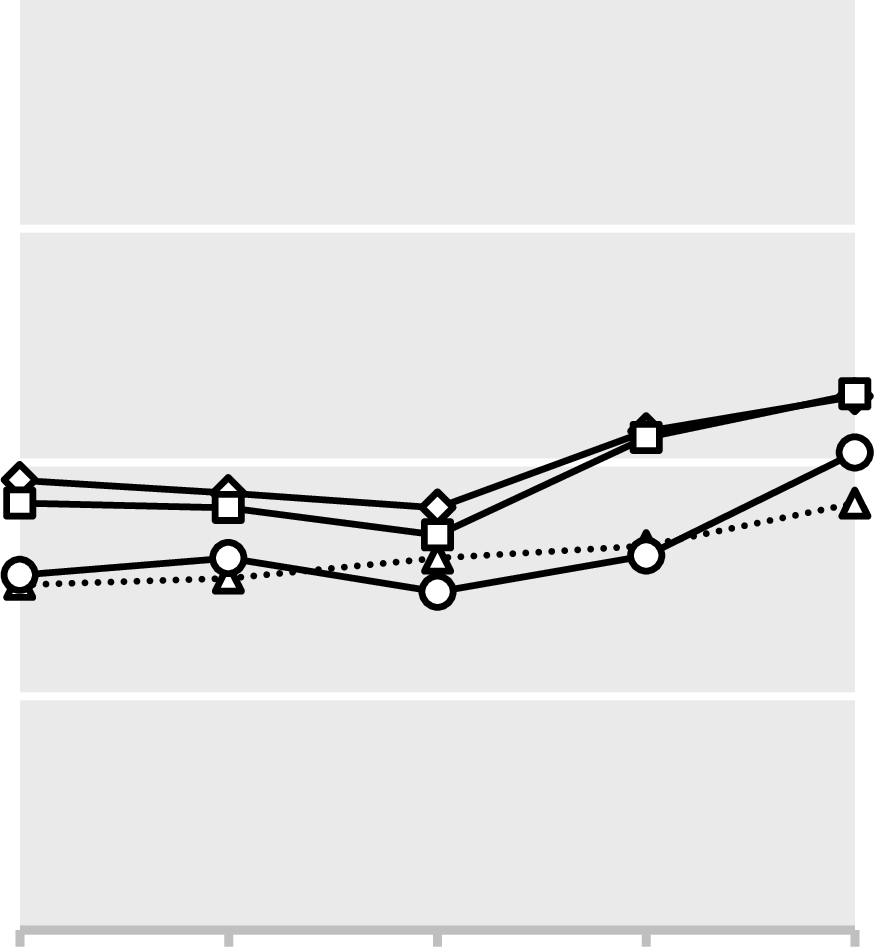}}
       \rput[b](0.117\linewidth,0.08){\scriptsize \textbf{\sf Network number}}
       \rput[bl](-0.05,4.8){\mbox{\scriptsize \textbf{{\sf Communication steps}}}}
       \rput[r](0.3,0.79){\scriptsize $\mathsf{10^0}$}
       \rput[r](0.3,1.73){\scriptsize $\mathsf{10^1}$}
       \rput[r](0.3,2.66){\scriptsize $\mathsf{10^2}$}
       \rput[r](0.3,3.61){\scriptsize $\mathsf{10^3}$}
       \rput[r](0.3,4.55){\scriptsize $\mathsf{10^4}$}

       \rput[t](0.429,0.58){\scriptsize $\mathsf{1}$}
       \rput[t](1.28,0.58){\scriptsize $\mathsf{2}$}
       \rput[t](2.132,0.58){\scriptsize $\mathsf{3}$}
       \rput[t](2.98,0.58){\scriptsize $\mathsf{4}$}
       \rput[t](3.825,0.58){\scriptsize $\mathsf{5}$}

       \rput[lt](1.7,2.01){\scriptsize \textbf{\sf D-ADMM}}
       \rput[t](3.6,2.37){\scriptsize \textbf{\sf \cite{Oreshkin10-OptimizationAnalysisDistrAveraging}}}
       \rput[l](1.1,3.0){\scriptsize \textbf{\sf \cite{Zhu09-DistributedInNetworkChannelCoding}}}
       \psline[linewidth=0.5pt](1.09,2.86)(0.76,2.5)
       \rput[lb](0.42,2.7){\scriptsize \textbf{\sf \cite{Schizas08-ConsensusAdHocWSNsPartI}}}
       %\rput[rb](3.59,4.01){\scriptsize \textbf{\sf \cite{Olshevsky11-ConvergenceSpeedDistributedConsensusAveraging}}}
       %\psgrid
     \end{pspicture}
     }
     \subfigure[LASSO]{\label{SubFig:Lasso}
     \begin{pspicture}(0.23\linewidth,5.1cm)
       \rput(0.117\linewidth,2.6){\includegraphics[scale=0.24]{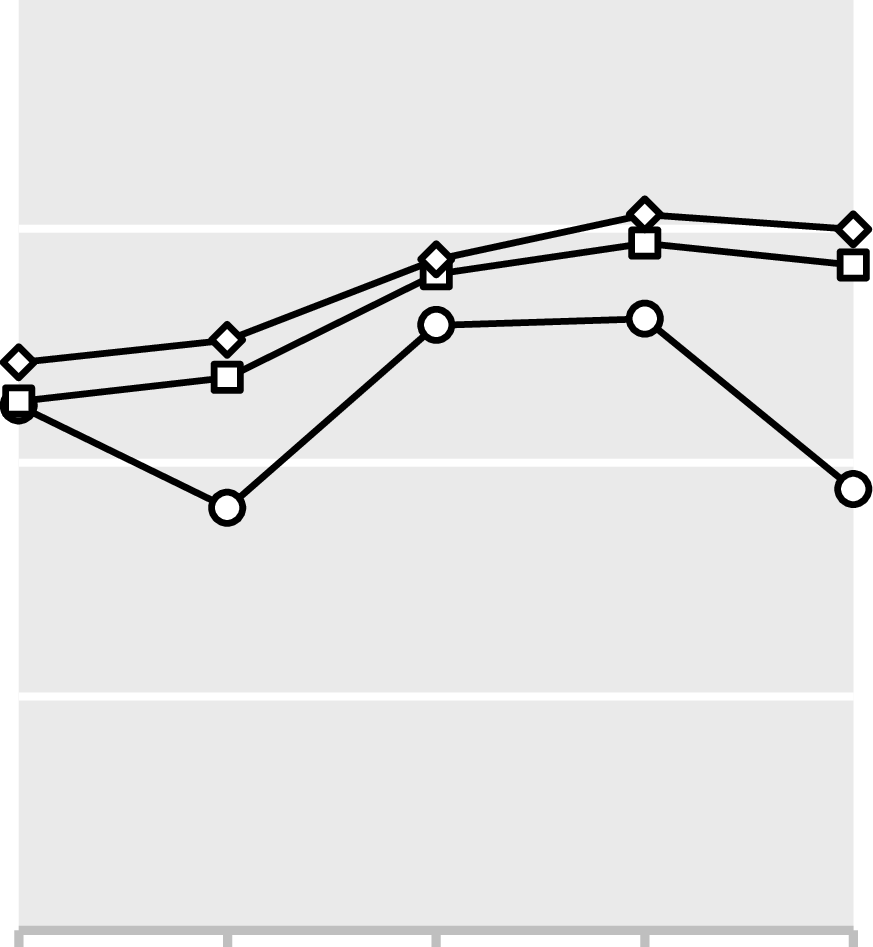}}
       \rput[b](0.117\linewidth,0.08){\scriptsize \textbf{\sf Network number}}
       \rput[bl](-0.05,4.8){\mbox{\scriptsize \textbf{{\sf Communication steps}}}}
       \rput[r](0.3,0.79){\scriptsize $\mathsf{10^0}$}
       \rput[r](0.3,1.73){\scriptsize $\mathsf{10^1}$}
       \rput[r](0.3,2.66){\scriptsize $\mathsf{10^2}$}
       \rput[r](0.3,3.61){\scriptsize $\mathsf{10^3}$}
       \rput[r](0.3,4.55){\scriptsize $\mathsf{10^4}$}

       \rput[t](0.429,0.58){\scriptsize $\mathsf{1}$}
       \rput[t](1.28,0.58){\scriptsize $\mathsf{2}$}
       \rput[t](2.132,0.58){\scriptsize $\mathsf{3}$}
       \rput[t](2.98,0.58){\scriptsize $\mathsf{4}$}
       \rput[t](3.825,0.58){\scriptsize $\mathsf{5}$}

       \rput[rt](3.6,2.55){\scriptsize \textbf{\sf D-ADMM}}
       \rput[lb](0.43,3.15){\scriptsize \textbf{\sf \cite{Schizas08-ConsensusAdHocWSNsPartI}}}
       \rput[l](1.255,3.58){\scriptsize \textbf{\sf \cite{Zhu09-DistributedInNetworkChannelCoding}}}
       \psline[linewidth=0.5pt](1.26,3.43)(1.0,3.0)
       %\psgrid
     \end{pspicture}
     }
     \subfigure[BPDN]{\label{SubFig:BPDN}
     \begin{pspicture}(0.23\linewidth,5.1cm)
       \rput(0.117\linewidth,2.6){\includegraphics[scale=0.24]{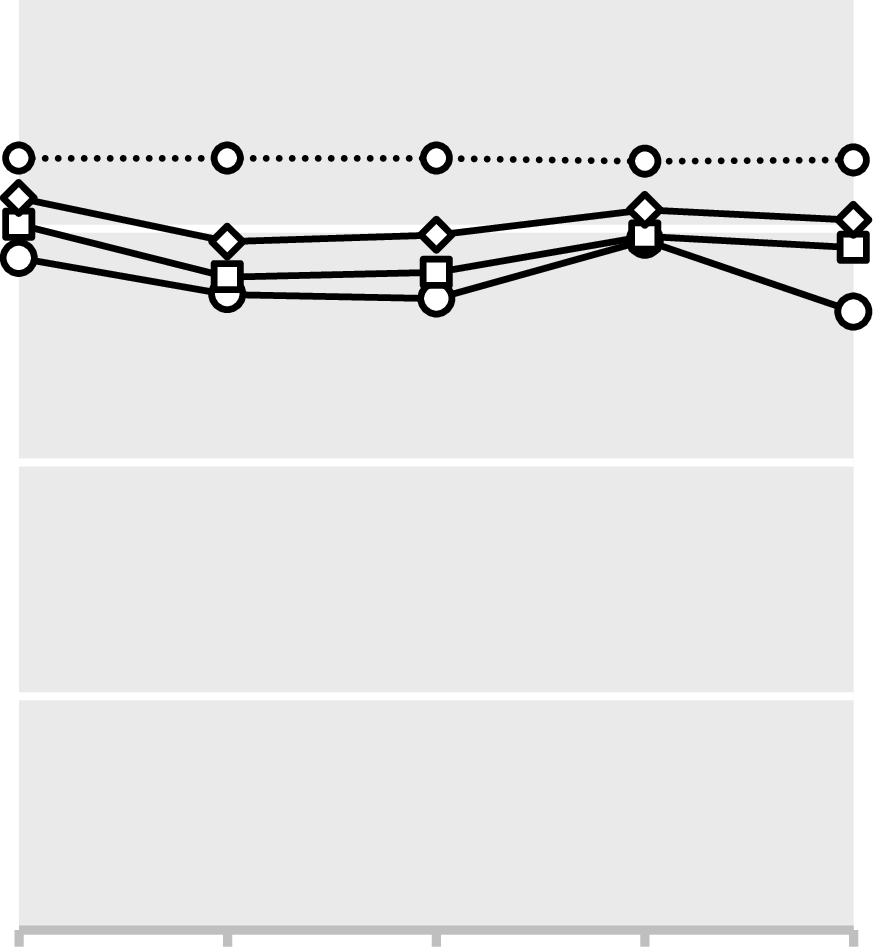}}
       \rput[b](0.117\linewidth,0.08){\scriptsize \textbf{\sf Network number}}
       \rput[bl](-0.05,4.8){\mbox{\scriptsize \textbf{{\sf Communication steps}}}}
       \rput[r](0.3,0.79){\scriptsize $\mathsf{10^0}$}
       \rput[r](0.3,1.73){\scriptsize $\mathsf{10^1}$}
       \rput[r](0.3,2.66){\scriptsize $\mathsf{10^2}$}
       \rput[r](0.3,3.61){\scriptsize $\mathsf{10^3}$}
       \rput[r](0.3,4.55){\scriptsize $\mathsf{10^4}$}

       \rput[t](0.429,0.58){\scriptsize $\mathsf{1}$}
       \rput[t](1.28,0.58){\scriptsize $\mathsf{2}$}
       \rput[t](2.132,0.58){\scriptsize $\mathsf{3}$}
       \rput[t](2.98,0.58){\scriptsize $\mathsf{4}$}
       \rput[t](3.825,0.58){\scriptsize $\mathsf{5}$}

       \rput[rt](3.76,3.2){\scriptsize \textbf{\sf D-ADMM}}
       \rput[lt](1.1,3.1){\scriptsize \textbf{\sf \cite{Schizas08-ConsensusAdHocWSNsPartI}}}
       \psline[linewidth=0.5pt](1.1,3.12)(0.88,3.6)
			 \rput[rt](0.85,3.1){\scriptsize \textbf{\sf \cite{Zhu09-DistributedInNetworkChannelCoding}}}
       \psline[linewidth=0.5pt](0.65,3.12)(0.7,3.52)
			 \rput[rb](3.77,3.97){\scriptsize \textbf{\sf \cite{Mateos10-DistributedSparseLinearRegression}}}
       %\psgrid
     \end{pspicture}
     }
     \subfigure[SVM]{\label{SubFig:SVM}
     \begin{pspicture}(0.23\linewidth,5.1cm)
       \rput(0.117\linewidth,2.6){\includegraphics[scale=0.24]{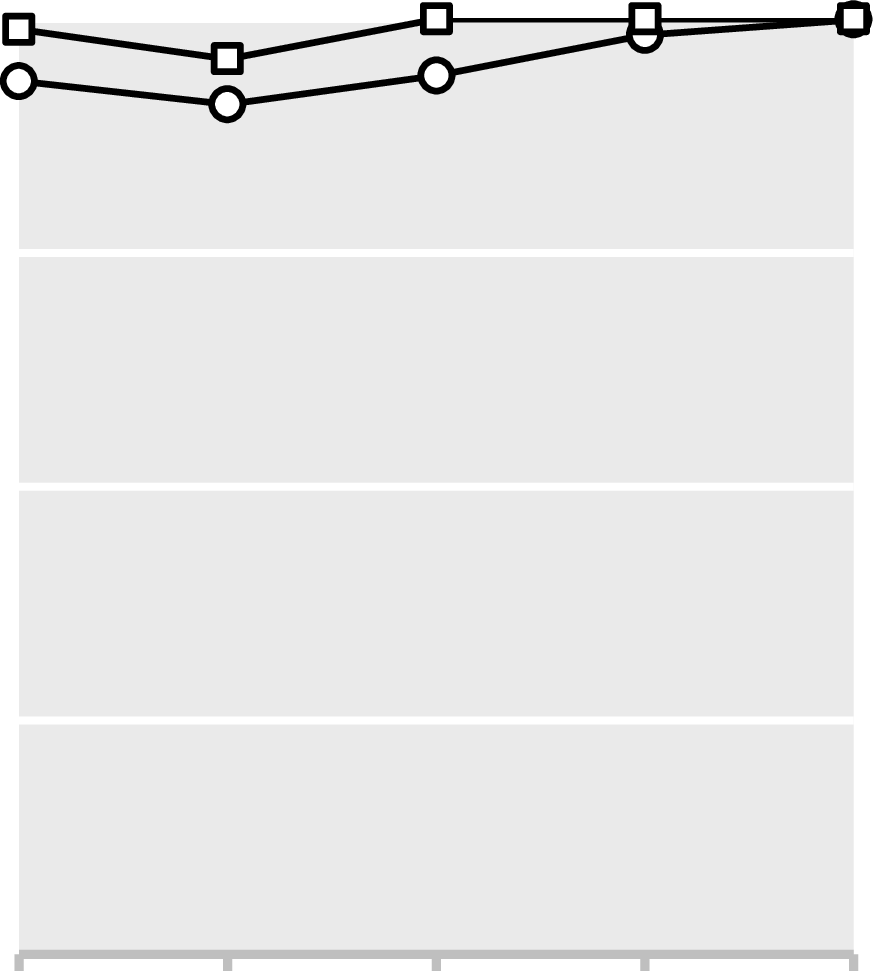}}

       \rput[b](0.117\linewidth,0.08){\scriptsize \textbf{\sf Network number}}
       \rput[bl](-0.05,4.8){\mbox{\scriptsize \textbf{{\sf Communication steps}}}}
       \rput[r](0.3,0.79){\scriptsize $\mathsf{10^0}$}
       \rput[r](0.3,1.73){\scriptsize $\mathsf{10^1}$}
       \rput[r](0.3,2.66){\scriptsize $\mathsf{10^2}$}
       \rput[r](0.3,3.61){\scriptsize $\mathsf{10^3}$}
       \rput[r](0.3,4.55){\scriptsize $\mathsf{10^4}$}

       \rput[t](0.429,0.58){\scriptsize $\mathsf{1}$}
       \rput[t](1.28,0.58){\scriptsize $\mathsf{2}$}
       \rput[t](2.132,0.58){\scriptsize $\mathsf{3}$}
       \rput[t](2.98,0.58){\scriptsize $\mathsf{4}$}
       \rput[t](3.825,0.58){\scriptsize $\mathsf{5}$}

       \rput[lt](0.67,4.04){\scriptsize \textbf{\sf D-ADMM}}
       \rput[lt](1.92,4.1){\scriptsize \textbf{\sf \cite{Zhu09-DistributedInNetworkChannelCoding}}}
       \psline[linewidth=0.5pt](1.92,4.12)(1.7,4.37)
       %\psgrid
     \end{pspicture}
     }
     \vspace{-0.1cm}
     \caption{Results of the simulations for $\text{(a)}$ consensus, $\text{(b)}$ LASSO, $\text{(c)}$ BPDN, and $\text{(d)}$ SVM.}
     \label{Fig:Experiments}
	\end{figure*}
  }

	\mypar{Performance measure: communication steps} We say that a \textit{Communication Step} (CS) has occurred when all the nodes have transmitted a vector of size~$n$ to its neighbors. All the algorithms we consider here, including D-ADMM, consist of one iterative loop. One iteration of D-ADMM, as well as of~\cite{Zhu09-DistributedInNetworkChannelCoding}, corresponds to one CS; one iteration of~\cite{Schizas08-ConsensusAdHocWSNsPartI} corresponds to two CSs, since each node transmits two vectors of size~$n$ per iteration. The number of CSs is proportional to the number of total communications. Thus, in a wireless scenario, the smaller the number of CSs, the lower the energy consumption.

	Note, however, that the CS measure does not take into account the computational complexity at each node. (Actually, D-ADMM, and the algorithms in~\cite{Schizas08-ConsensusAdHocWSNsPartI,Zhu09-DistributedInNetworkChannelCoding} have similar computational complexities.) Also, this measure is not necessarily related with execution time. In fact, while D-ADMM requires less CSs than competing algorithms (as we will see), it may be slower than some of them. The reason is because D-ADMM is asynchronous, while all the other algorithms are synchronous. Scenarios allowing synchronous transmissions are, however, limited to very controlled environments, such as computer clusters or super-computers, where approaches like~\cite[\S7.2]{Boyd11-ADMM} would probably be more appropriate than distributed algorithms. On the other hand, in wireless networks, the single fact that one node cannot transmit and receive messages at the same time forces synchronous algorithms to operate asynchronously.
	%Execution time depends strongly on the platform,  communication protocol, and implementation (e.g., programming language), among other factors. We will focus instead on CSs, which are intrinsic to the algorithms and independent of these factors.

	\isdraft{                                                                 % DRAFT
    \begin{table}                                                           % DRAFT
      %\vspace{-1cm}
      \centering                                                            % DRAFT
      \caption{Network models}                                              % DRAFT
      \vspace{-0.3cm}                                                       % DRAFT
      \label{Tab:NetworkModels}                                             % DRAFT
      \includegraphics[scale=0.9]{figures/SingleCol/NetworksTable.eps}      % DRAFT
      \vspace{-0.5cm}                                                       % DRAFT
    \end{table}                                                             % DRAFT
  }{
	\begin{table}
    \centering
    \caption{Network models}
    \footnotesize
    \label{Tab:NetworkModels}
    \renewcommand{\arraystretch}{1.2}
    %\resizebox{0.7\linewidth}{!}{
    \begin{tabular}{@{}clllc@{}}
    \toprule[1pt]
    Network && Model (parameters) && $\#$ Colors \\
    \midrule

    1 && Erd\H os-R\'enyi \,\,($0.12$) && \phantom{0}5\\
    2 && Watts-Strogatz \,\,($4,0.4$) && \phantom{0}4 \\
    3 && Barabasi \,\,($2$) && \phantom{0}3\\
    4 && Geometric \,\,($0.23$) && 10\\
    5 && Lattice \,\,($5\times 10$) && \phantom{0}2\\
    \bottomrule[1pt]
    \end{tabular}
    %}
  \end{table}
  }

% 	\begin{table*}
%     \centering
%     \renewcommand{\arraystretch}{1.2}
%     %\resizebox{0.7\linewidth}{!}{
%     \begin{tabular}{@{}llrcccc@{}}
%     \toprule[1pt]
%     Application && Max CSs && $\epsilon$ && $\rho$ \\
%     \midrule
%     Consensus && $10^3$ && $10^{-4}$         && $\{10^{-4}, 10^{-3} 10^{-2} 10^{-1}, 1, 10, 10^{2}\}$ \\
%     Lasso     && $10^3$ && $5\times 10^{-3}$ && $\{10^{-3}, 10^{-2} 10^{-1}, 1, 10\}$ \\
%     BPDN      && $2\times10^3$ && $10^{-4}$  && $\{10^{-4}, 10^{-3} 10^{-2} 10^{-1}, 1, 10, 10^{2}\}$ \\
%     SVM       && $10^3$ && $10^{-3}$         && $\{10^{-4}, 10^{-3} 10^{-2} 10^{-1}, 1, 10, 10^{2}\}$ \\
%     \bottomrule[1pt]
%     \end{tabular}
%     %}
%    \end{table*}
%

	\mypar{Experimental setup}
	We generated~$5$ networks with~$P = 50$ nodes according to the models of Table~\ref{Tab:NetworkModels}; see~\cite{Mota12-DistributedBP} for a description of these models. Table~\ref{Tab:NetworkModels} also gives the number of colors for each network. The results of our simulations are in Fig.~\ref{Fig:Experiments}, where each plot depicts the number of CSs as a function of the network. Having computed the solution~$x^\star$ of~\eqref{Eq:IntroSeparableProb} beforehand and in a centralized way, each algorithm stopped whenever $\|x^k - x^\star\|/\|x^\star\| \leq \epsilon$, or when a maximum number of~$M$ CSs were reached. In the case of consensus, BPDN, and SVM, $x^k$ denotes the estimate of~$x^\star$ at an arbitrary node; in the case of LASSO, it represents the global estimate of the network, since each node only estimates some components of~$x^\star$. We used the following values for the pair~$(\epsilon,M)$: $(10^{-4},10^3)$ for consensus, $(5\times 10^{-3}, 10^3)$ for LASSO, $(10^{-4},2\times 10^3)$ for BPDN, and $(10^{-3},10^3)$ for SVM. Since the problem in step~\ref{SubAlg:ADMMGenProb} of Algorithm~\ref{Alg:DADMMGeneral} does not have a closed-form solution for all the applications we consider, except for consensus, it has to be solved iteratively in each of the algorithms we compare. To make a fair comparison in terms of CSs, we thus use the same solver in all the algorithms, i.e., the problem in step~\ref{SubAlg:ADMMGenProb} of Algorithm~\ref{Alg:DADMMGeneral} is solved with the same precision in all the algorithms we compare.

	It is known that the parameter~$\rho$ affects strongly the performance of ADMM-based algorithms. Hence, to make a fair comparison, we ran each (ADMM-based) algorithm for several values of~$\rho$ and chose always the best result, i.e., the smaller number of CSs. The values for~$\rho$ were taken from the set $\{10^{-4}, 10^{-3} 10^{-2} 10^{-1}, 1, 10, 10^{2}\}$.

	\mypar{Consensus} For the consensus problem, we generated each~$\theta_p$ randomly from a Gaussian distribution: $\theta_p \overset{\text{i.i.d.}}{\sim} \mathcal{N}(10,10^4)$. Fig.~\ref{SubFig:Consensus} shows the results for D-ADMM, the ADMM-based algorithms~\cite{Schizas08-ConsensusAdHocWSNsPartI,Zhu09-DistributedInNetworkChannelCoding}, and the algorithm~\cite{Oreshkin10-OptimizationAnalysisDistrAveraging}, which is considered to be the fastest consensus algorithm~\cite{Erseghe11-FastConsensusByADMM}. Note that~\cite{Oreshkin10-OptimizationAnalysisDistrAveraging}  was designed for consensus only and cannot be easily generalized to solve~\eqref{Eq:IntroSeparableProb}. Fig.~\ref{SubFig:Consensus} shows that D-ADMM has a performance very similar to that of~\cite{Oreshkin10-OptimizationAnalysisDistrAveraging}.

	\mypar{LASSO and BPDN} The matrix~$A$ for the problems LASSO and BPDN was taken from problem~$902$ of the Sparco toolbox~\cite{Sparco}. The vector~$b$ was generated as~$b = As + n$, where~$s$ is a sparse vector and~$n$ is Gaussian noise. We chose~$\sigma = 0.1$ and $\beta = 0.3$ for the noise parameters, and~$\delta = 10^{-2}$ for the approximation parameter in LASSO. The results of these experiments for LASSO and BPDN are shown, respectively, in Figs.~\ref{SubFig:Lasso} and~\ref{SubFig:BPDN}. Additionally, we show the performance of Algorithm~$3$ of~\cite{Mateos10-DistributedSparseLinearRegression}, which is an ADMM-based algorithm specifically designed to solve BPDN. This algorithm has the advantage of requiring much simpler computations at each node, but in our simulations it achieved the maximum number of CSs in all but the last two networks. In both LASSO and BPDN, D-ADMM was always the algorithm requiring fewer CSs to converge.

	\mypar{SVM} For the SVM problem~\eqref{Eq:SVM}, we used data from~\cite{UCIMachineLearningRepository}, namely two overlapping sets of points from the Iris dataset. The parameter~$\beta$ was always set to~$1$. In this case, the algorithm from~\cite{Schizas08-ConsensusAdHocWSNsPartI} achieved always the maximum number of CSs and thus is not represented in Fig.~\ref{SubFig:SVM}, which shows the simulation results. Again, we see that D-ADMM was the algorithm requiring the smallest number of CSs to converge.

\section{Conclusions}
\label{Sec:Conclusions}

  We proposed an algorithm for solving separable problems in networks, in a distributed way. Each node has a private cost and a private constraint set, but all nodes cooperate to solve the optimization problem that minimizes the sum of all costs and that has the intersection of all sets as a constraint. Our algorithm hinges on a coloring scheme of the network, according to which the nodes operate asynchronously. This results in an algorithm with fewer communication requirements than previous algorithms, as shown experimentally for several problems. Although we proved the convergence of the algorithm, it still remains an open question to explain theoretically why the algorithm is more efficient than previous algorithms.

\bibliographystyle{IEEEbib}

{ \isdraft{\singlespace}{}
\bibliography{paper}

\begin{thebibliography}{10}

\bibitem{Mota12-DistributedBP}
J.~{Mota}, J.~{Xavier}, P.~{Aguiar}, and M.~{P\"uschel},
\newblock ``Distributed basis pursuit,''
\newblock {\em IEEE Trans. Sig. Proc.}, vol. 60, no. 4, 2012.

\bibitem{Akyildiz02-WirelessNetworksASurvey}
I.~{Akyildiz}, Y.~{Sankarasubramaniam}, and E.~{Cayirci},
\newblock ``Wireless sensor networks: a survey,''
\newblock {\em Comp. Netw.}, vol. 38, 2002.

\bibitem{Fischione11-DesignPrinciplesWSN}
P.~{Fischione}, P.~{Park}, and K.~{Johansson},
\newblock {\em Wireless Network Based Control}, chapter Design Principles of
  Wireless Sensor Network Protocols for Control Applications,
\newblock Springer, 2011.

\bibitem{Boyd11-ADMM}
S.~{Boyd}, N.~{Parikh}, E.~{Chu}, B.~{Peleato}, and J.~{Eckstein},
\newblock ``Distributed optimization and statistical learning via the
  alternating method of multipliers,''
\newblock {\em Found. Trends Mach. Learn.}, vol. 3, no. 1, pp. 1--122, 2011.

\bibitem{Han12-NoteOnADMM}
D.~{Han} and X.~{Yuan},
\newblock ``A note on the alternating direction method of multipliers,''
\newblock {\em J. Optim. Theory Appl.}, 2012.

\bibitem{Tsitsiklis86-DistrAsynchronous}
J.~{Tsitsiklis}, D.~{Bertsekas}, and M.~{Athans},
\newblock ``Distributed asynchronous deterministic and stochastic gradient
  optimization algorithms,''
\newblock {\em IEEE Trans. Aut. Contr.}, vol. AC-31, no. 9, 1986.

\bibitem{Rabbat04-DistributedOptimizationSensorNetworks}
M.~{Rabbat} and R.~{Nowak},
\newblock ``Distributed optimization in sensor networks,''
\newblock in {\em Proc. IPSN'04}, 2004, pp. 20--27.

\bibitem{CXVOptimInSPandComm}
P~{Palomar} and Y.~{Eldar},
\newblock {\em Convex Optimization in Signal Processing and Communications},
\newblock Cambridge Univ. Press, 2010.

\bibitem{Bertsekas97-ParallelDistributed}
D.~{Bertsekas} and J.~{Tsitsiklis},
\newblock {\em Parallel and Distributed Computation: Numerical Methods},
\newblock Athena Scientific, 1997.

\bibitem{Jakovetic11-CooperativeConvexOptimization}
D.~{Jakoveti\'c}, J.~{Xavier}, and J.~{Moura},
\newblock ``Cooperative convex optimization in networked systems: Augmented
  lagrangian algorithms with directed gossip communication,''
\newblock {\em IEEE Trans. Sig. Proc.}, vol. 59, no. 8, 2011.

\bibitem{Rabbat05-GeneralizedConsensusAlgs}
M.~{Rabbat}, R.~{Nowak}, and J.~{Bucklew},
\newblock ``Generalized consensus algorithms in networked systems with erasure
  links,''
\newblock in {\em IEEE Workshop Sig. Proc. Advances Wirel. Comm.}, 2005.

\bibitem{Schizas08-ConsensusAdHocWSNsPartI}
I.~{Schizas}, A.~{Ribeiro}, and G.~{Giannakis},
\newblock ``Consensus in \textit{ad hoc} wsns with noisy links - \text{Part I}:
  Distributed estimation of deterministic signals,''
\newblock {\em IEEE Trans. Sig. Proc.}, vol. 56, no. 1, pp. 350--364, 2008.

\bibitem{Zhu09-DistributedInNetworkChannelCoding}
H.~{Zhu}, G.~{Giannakis}, and A.~{Cano},
\newblock ``Distributed in-network channel decoding,''
\newblock {\em IEEE Trans. Sig. Proc.}, vol. 57, no. 10, 2009.

\bibitem{Bazerque10-DistributedSpectrumSensing}
J.~{Bazerque} and G.~{Giannakis},
\newblock ``Distributed spectrum sensing for cognitive radio networks by
  exploiting sparsity,''
\newblock {\em IEEE Trans. Sig. Proc.}, vol. 58, no. 3, 2010.

\bibitem{Forero10-ConsensusBasedDistributedSVMs}
P.~{Forero}, A.~{Cano}, and G.~{Giannakis},
\newblock ``Consensus-based distributed support vector machines,''
\newblock {\em J.M.L.R.}, vol. 11, 2010.

\bibitem{Mateos10-DistributedSparseLinearRegression}
G.~{Mateos}, J.~{Bazerque}, and G.~{Giannakis},
\newblock ``Distributed sparse linear regression,''
\newblock {\em IEEE Trans. Sig. Proc.}, vol. 58, no. 10, 2010.

\bibitem{Bazerque11-GroupLasso}
J.~{Bazerque}, G.~{Mateos}, and G.~{Giannakis},
\newblock ``Group-\text{L}asso on splines for spectrum cartography,''
\newblock {\em IEEE Trans. Sig. Proc.}, vol. 59, no. 10, 2011.

\bibitem{Erseghe11-FastConsensusByADMM}
T.~{Erseghe}, D.~{Zennaro}, E.~{Dall'Anese}, and L.~{Vangelista},
\newblock ``Fast consensus by the alternating direction multipliers method,''
\newblock {\em IEEE Trans. Sig. Proc.}, vol. 59, no. 11, 2011.

\bibitem{Mota11-ADMMProof}
J.~{Mota}, J.~{Xavier}, P.~{Aguiar}, and M.~{P\"uschel},
\newblock ``A proof of convergence for the alternating direction method of
  multipliers applied to polyhedral-constrained functions,''
\newblock \url{http://arxiv.org/abs/1112.2295}, 2011.

\bibitem{He12-ADMGaussianBack}
B.~{He}, M.~{Tao}, and X.~{Yuan},
\newblock ``Alternating direction method with $\textrm{Gaussian}$ back
  substitution for separable convex programmming,''
\newblock {\em SIAM J. Optim.}, vol. 22, no. 2, 2012.

\bibitem{Luo12-LinearConvergenceADMM}
Z.~{Luo},
\newblock ``On the linear convergence of the alternating direction method of
  multipliers,''
\newblock \url{arxiv.org/abs/1208.3922}, 2012.

\bibitem{Oreshkin10-OptimizationAnalysisDistrAveraging}
B.~{Oreshkin}, M.~{Coates}, and M.~{Rabbat},
\newblock ``Optimization and analysis of distributed averaging with short node
  memory,''
\newblock {\em IEEE Trans. Sig. Proc.}, vol. 58, no. 5, 2010.

\bibitem{Tibshirani96-RegressionShrinkageLasso}
R.~{Tibshirani},
\newblock ``Regression shrinkage and selection via the lasso,''
\newblock {\em J. R. Statist. Soc. B}, vol. 58, no. 1, 1996.

\bibitem{AtomicDecompBP}
S.~{Chen}, D.~{Donoho}, and M.~{Saunders},
\newblock ``Atomic decomposition by basis pursuit,''
\newblock {\em SIAM J. Sc. Cp.}, vol. 20, no. 1, 1998.

\bibitem{Friedlander07ExactRegularizationConvexPrograms}
M.~{Friedlander} and P.~{Tseng},
\newblock ``Exact regulatization of convex programs,''
\newblock {\em SIAM J. Optim.}, vol. 18, no. 4, 2007.

\bibitem{Bishop06-PatternRecognitionMachineLearning}
C.~{Bishop},
\newblock {\em Pattern Recognition and Machine Learning},
\newblock Springer, 2006.

\bibitem{Sparco}
E.~{Berg}, M.~{Friedlander}, G.~{Hennenfent}, F.~{Herrmann}, R.~{Saab}, and
  \"O. {Yilmaz},
\newblock ``Sparco: a testing framework for sparse reconstruction,''
\newblock Tech. {R}ep., Dept. Computer Science, University of British Columbia,
  Vancouver, 2007.

\bibitem{UCIMachineLearningRepository}
A.~{Frank} and A.~{Asuncion},
\newblock {\em UCI Machine Learning Repository},
\newblock Univ. Calif., 2010.

\end{thebibliography}
}

\begin{IEEEbiography}
[{\includegraphics[width=1in,clip]{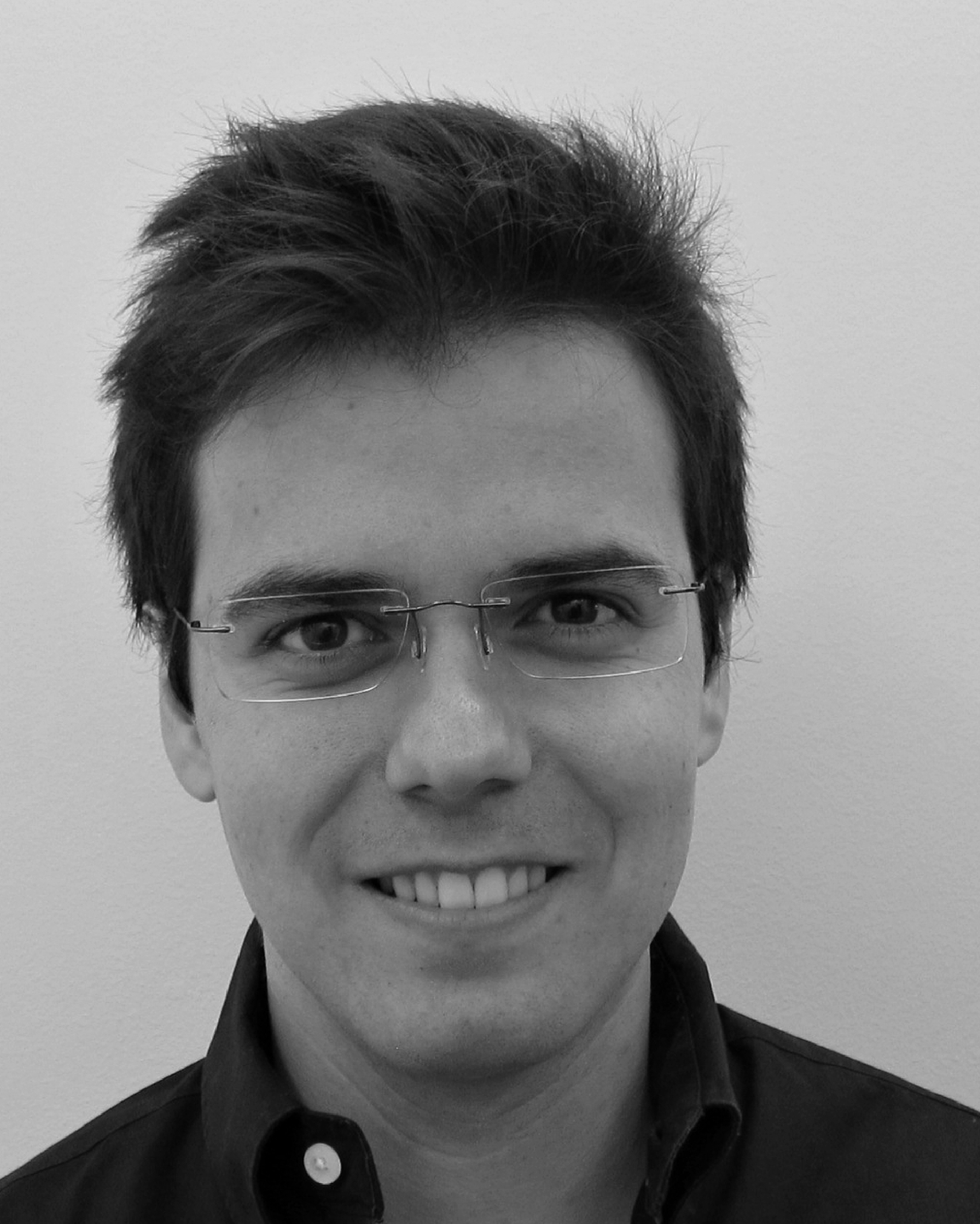}}] {Jo\~ao Mota} received a M.S. degree in Electrical and Computer Engineering from Instituto Superior T\'ecnico, Technical University of Lisbon, Lisbon, Portugal, in 2008. He is currently working towards his Ph.D. degree in Electrical and Computer Engineering, within a joint program between Carnegie Mellon University, Pittsburgh, PA, and Instituto Superior T\'ecnico, Lisbon, Portugal. His research interests include distributed optimization and control, and sensor networks.
\end{IEEEbiography}

\begin{IEEEbiography}
[{\includegraphics[width=1in,height=1.25in,clip,keepaspectratio]%
  {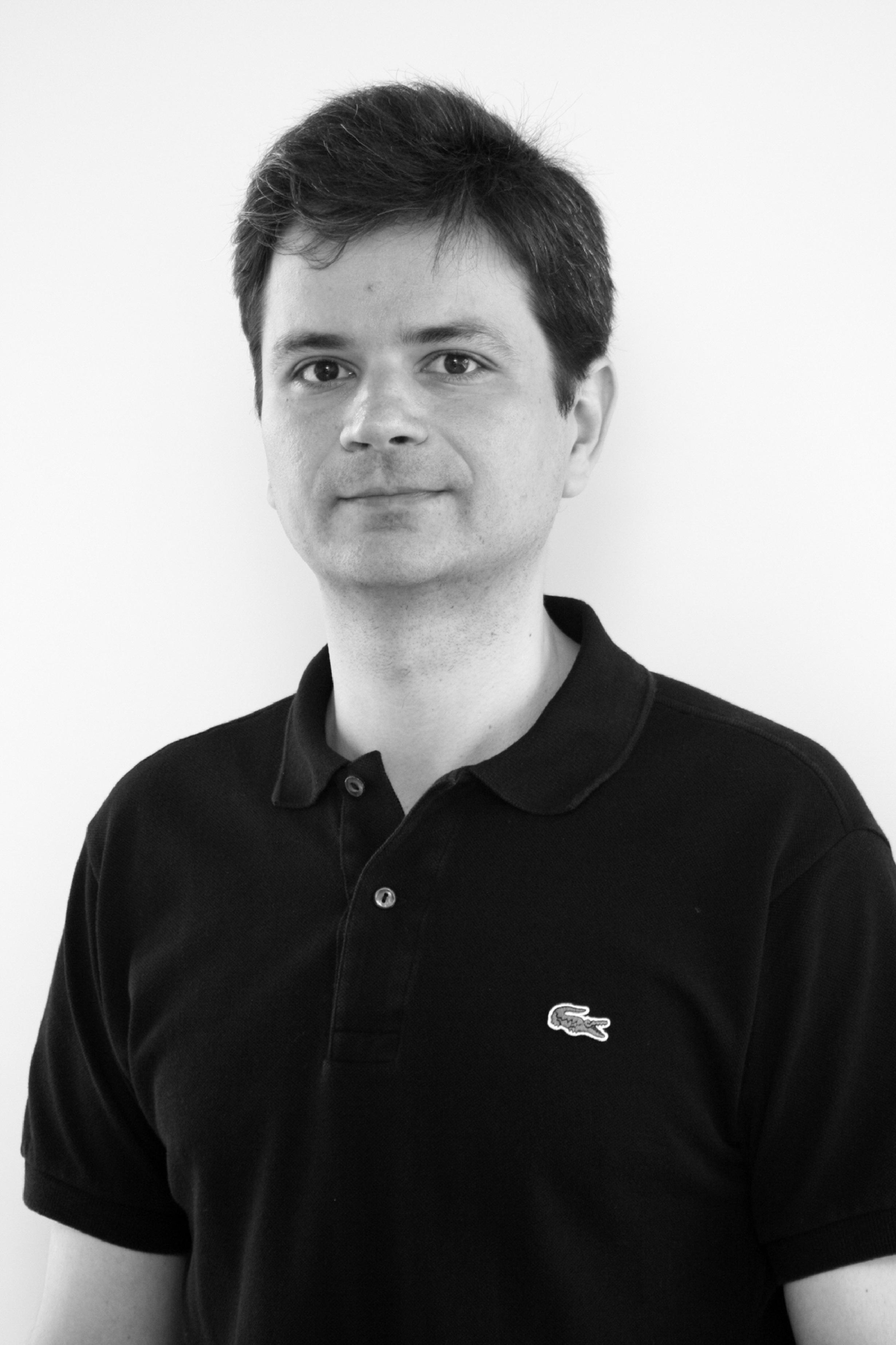}}] {Jo\~ao Xavier} (S'97--M'03) received the Ph.D. degree in Electrical and Computer Engineering from Instituto Superior Tecnico (IST), Lisbon, Portugal, in 2002. Currently, he is an Assistant Professor in the Department of Electrical and Computer Engineering, IST. He is also a Researcher at the Institute of Systems and Robotics (ISR), Lisbon, Portugal. His current research interests are in the area of optimization, sensor networks and signal processing on manifolds.
\end{IEEEbiography}

\begin{IEEEbiography}
[{\includegraphics[width=1in,height=1.25in,clip,keepaspectratio]{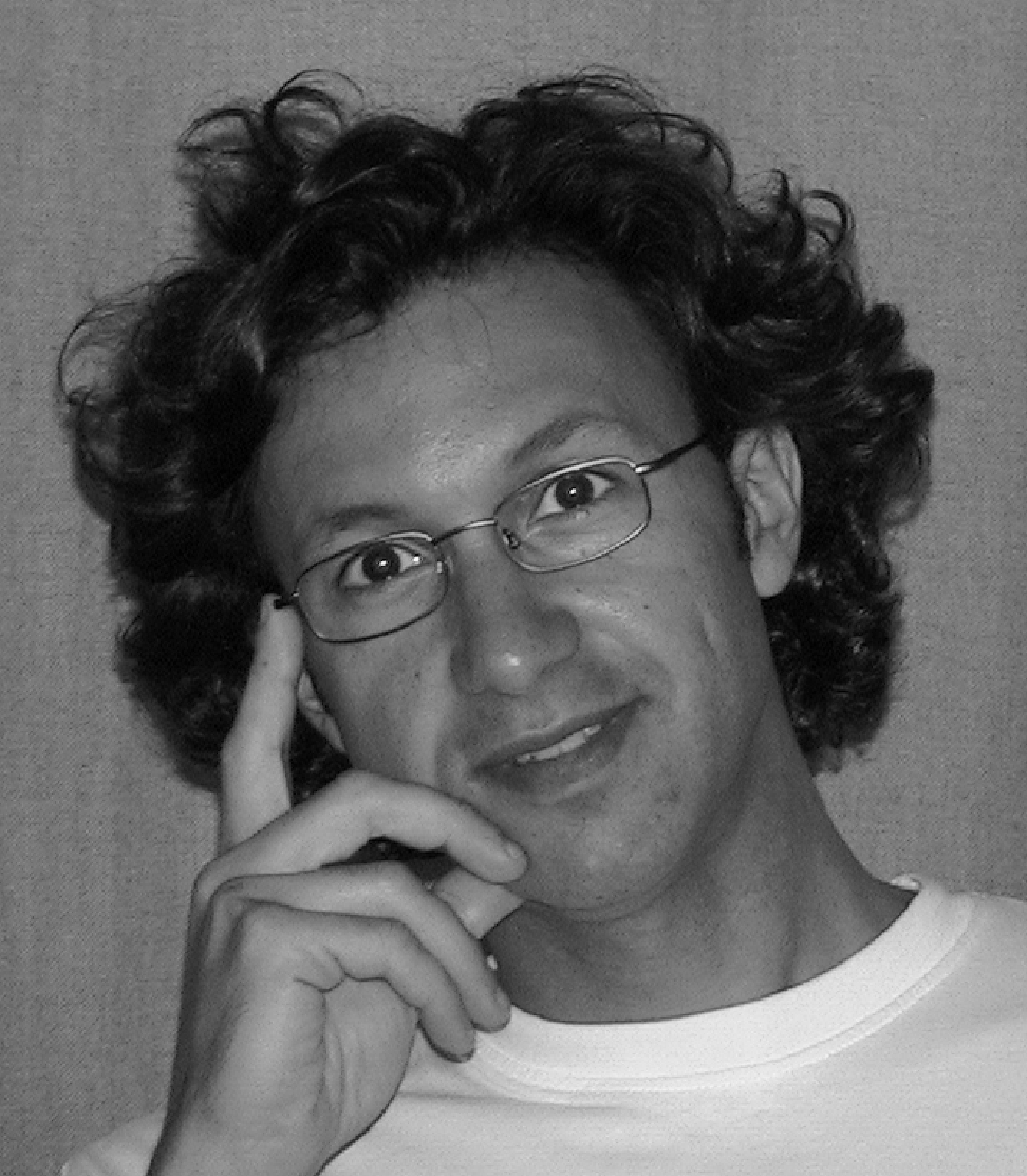}}]
{Pedro M.~Q.~Aguiar} (S'95--M'00--SM'08) received the Ph.D. degree in
electrical and computer engineering from the Instituto Superior
T\'ecnico, Technical University of Lisbon, Lisbon, Portugal, in 2000.

He is currently an Assistant Professor with the Instituto Superior
Técnico, Technical University of Lisbon, Lisbon, Portugal. He is also
affiliated with the Institute for Systems and Robotics, Lisbon,
Portugal, and has been Visiting Scholar with Carnegie-Mellon
University, Pittsburgh, PA, and a Consultant with Xerox Palo Alto
Research Center, Palo Alto, CA. His main research interests are in
image analysis and computer vision.
\end{IEEEbiography}

\begin{IEEEbiography}
[{\includegraphics[width=1in,height=1.25in,clip,keepaspectratio]%
  {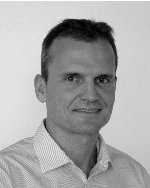}}] {Markus P{\"u}schel} (M'99--SM'05)
  is a Professor of Computer Science at ETH Zurich, Switzerland.
  Before, he was a Professor of Electrical and Computer Engineering at
  Carnegie Mellon University, where he still has an adjunct status.
  He received his Diploma (M.Sc.) in Mathematics and his Doctorate (Ph.D.)
  in Computer Science, in 1995 and 1998, respectively, both from the
  University of Karlsruhe, Germany. From 1998-1999 he was a Postdoctoral
  Researcher at Mathematics and Computer Science, Drexel University.
  From 2000-2010 he was with Carnegie Mellon University, and since 2010
  he has been with ETH Zurich. He was an Associate Editor for the IEEE
  Transactions on Signal Processing, the IEEE Signal Processing Letters,
  was a Guest Editor of the Proceedings of the IEEE and the Journal of
  Symbolic Computation, and served on various program committees of conferences
  in computing, compilers, and programming languages. He is a recipient of the
  Outstanding Research Award of the College of Engineering at Carnegie Mellon and
  the Eta Kappa Nu Award for Outstanding Teaching. He also holds the title of
  Privatdozent at the University of Technology, Vienna, Austria.
  In 2009 he cofounded SpiralGen,~Inc.
\end{IEEEbiography}

\end{document}